\documentclass[smallcondensed,envcountsect]{svjour3}
\usepackage{geometry}                
\usepackage{graphicx}
\usepackage{amssymb}
\usepackage{epstopdf}
\usepackage{algorithm2e}
\usepackage[reqno]{amsmath}
\usepackage{amssymb}
\usepackage[]{graphicx}               
\usepackage{bm}
\usepackage{esvect}
\usepackage{ulem}
\usepackage{adjustbox}
\usepackage{accents}
\usepackage{stmaryrd}
\usepackage{xcolor}
\usepackage{enumitem}
\usepackage{hyperref}
\usepackage{booktabs}
\usepackage{longtable}
\usepackage[section]{placeins}
\usepackage{tabu}
\usepackage{pgfplots}

\usepackage{color}
\definecolor{chcol}{rgb}{0.4,0.,0.9}

\newcommand{\halfComma}{\kern 0.083334em}

\newcommand{\avg}[1]{\left\{\hspace*{-3pt}\left\{#1\right\}\hspace*{-3pt}\right\}}
\newcommand{\jump}[1]{\ensuremath{\left\llbracket #1 \right\rrbracket}}
\newcommand{\dS}{{\,\operatorname{dS}}}         
\renewcommand{\sl}{+}                                      
\newcommand\iprod[1]{\left\langle #1\right\rangle}                                             
\newcommand\inorm[1]{\left |\left| #1\right|\right|}                                               
\newcommand\iprodN[1]{\left\langle #1\right\rangle_{\!N}}                                 
\newcommand\spacevec[1]{\accentset{\,\rightarrow}{#1}}                        
\newcommand\contravec[1]{\tilde{ #1}}                                                     
\newcommand\statevec[1]{\mathbf #1}                                                     
\newcommand\statevecGreek[1]{\boldsymbol #1}                                     
\newcommand\contrastatevec[1]{\tilde{\mathbf #1}}                                 
\newcommand\acclrvec[1]{\accentset{\,\leftrightarrow}{#1}}                      
\newcommand\bigstatevec[1]{\acclrvec{{\mathbf #1}}}                              
\newcommand\blockvec[1]{\acclrvec{{\mathbf #1}}}                                  
\newcommand\bigcontravec[1]{\acclrvec{\tilde{\mathbf{#1}}}}                   
\newcommand\vecNabla{\accentset{\,\rightarrow}\nabla}                         
\newcommand\vecNablaXi{\accentset{\,\rightarrow}\nabla_{\!\xi}}            
\newcommand\vecNablaX{\accentset{\,\rightarrow}\nabla_{\!x}}              
\newcommand\mmatrix[1]{\underbar{#1}}				
\newcommand\matrixvec[1]{\mathcal #1}                           
\newcommand\bigmatrix[1]{\mathfrak #1}                          
\newcommand{\dmat}{\matrixvec{D}}     

\newcommand\boundaryfaces{\genfrac{}{}{0pt}{}{\mathrm{boundary}}{\mathrm{faces}}}

\newcommand\DGfaces{\genfrac{}{}{0pt}{}{\mathrm{interior DG}}{\mathrm{faces}}}
\newcommand\sharedpoints{\genfrac{}{}{0pt}{}{\mathrm{shared}}{\mathrm{interface\; points}}}
\newcommand{\testfuncOne}{\statevecGreek{\varphi}}

\newcommand{\IN}[1]{\mathbb I^{N}\!\!\left(#1\right)} 
\newcommand{\PN}[1]{\mathbb P^{#1}}
\newcommand{\LTwo}[1]{\mathbb L^{2}\!\left(#1\right)}
\newcommand\oneHalf{\frac{1}{2}}

\usepackage{xargs}    
\usepackage[colorinlistoftodos,prependcaption,textsize=tiny]{todonotes}
\newcommandx{\unsure}[2][1=]{\todo[linecolor=blue,backgroundcolor=blue!25,bordercolor=blue,#1]{#2}}
\newcommandx{\changeThis}[2][1=]{\todo[linecolor=red,backgroundcolor=red!25,bordercolor=red,#1]{#2}}

\title{A Split-Form, Stable CG/DG-SEM for Wave Propagation Modeled by Linear Hyperbolic Systems}
\author{David A. Kopriva and Gregor J. Gassner}
\institute{David A. Kopriva \at  Professor Emeritus, Department of Mathematics, The Florida State University, Tallahassee, FL 32306, USA (\email{kopriva@math.fsu.edu}) and Computational Science Research Center, San Diego State University, San Diego, CA, USA\\
Gregor J. Gassner \at Department for Mathematics and Computer Science; Center for Data and Simulation Science,
University of Cologne, Weyertal 86-90, 50931, Cologne, Germany (\email{ggassner@uni-koeln.de})\\ 
           }

\begin{document}
\maketitle
\begin{abstract}
We present a hybrid continuous and discontinuous Galerkin spectral element approximation that leverages the advantages of each approach. The continuous Galerkin approximation is used on interior element faces where the equation properties are continuous. A discontinuous Galerkin approximation is used at physical boundaries and if there is a jump in properties at a face. The approximation uses a split form of the equations and two-point fluxes to ensure stability for unstructured quadrilateral/hexahedral meshes with curved elements. The approximation is also conservative and constant state preserving on such meshes. Spectral accuracy is obtained for all examples, which include wave scattering at a discontinuous medium boundary.
\end{abstract}
\section{Introduction}
Spectral element methods, which have desired features of low dissipation and dispersion errors \cite{hu1999}, \cite{Airnsworth2004}, \cite{gassner2011} and spectral or even exponential convergence rates \cite{CHQZ:2006}, have been used for many years to compute linear wave propagation problems. They have been used for electromagnetic wave \cite{Koprivaetal1999},\cite{Warburton1999},\cite{Hesthaven:2002uq},\cite{Cohen2006}, elastic 
wave \cite{wilcox2010}, linear acoustics \cite{Rasetarineraetal2001},\cite{ISI:000226090600009},\cite{PindEtAl2019}, and optics \cite{VanLith2017} problems, for instance. The discontinuous Galerkin Spectral Element Method  (DGSEM) version is robust, and in recent years provably stable forms have been developed for curved element meshes \cite{Gassner:2013ij}. One of the reasons for its success is the inter-element dissipation added at element faces introduced by the numerical fluxes. The DGSEM approximations are also naturally suitable for problems with material interfaces where properties jump.

Recently, a nodal Continuous (CG)/Discontinuous Galerkin (DG) hybrid method has been introduced and used in the electromagnetics community for time domain Maxwell's equations models \cite{Angulo:2015ja} under the name of CGDGTD (CG/DG Time Domain). The authors note that hybrid approach has a number of advantages. Parallelism, for instance, is enhanced  at processor boundaries over a pure CG approximation by exploiting the simplified coupling that DG elements have. On the other hand, in smooth regions and away from processor boundaries the CG approximation requires fewer degrees of freedom than DG, since nodes on both sides of a face do not have to be stored. The CG approximation also allows for a larger explicit time step, and has no inter-element dissipation. Combining the approaches, then, allows for the advantages of each to be exploited. Second and third order methods have been used with the hybrid approach.

The hybrid approach differs from a coupled CG/DG approximation, such as that proposed in \cite{DAWSON2004289}. In the mixed approach, some regions will be approximated by DG and other by CG, with a coupling condition when the two approximations meet. In \cite{DAWSON2004289}, for instance, the DG regions approximated a first order system, while in CG regions a second order equation was approximated. Instead, the hybrid approach approximates the first order system everywhere and the solutions and test functions are allowed to be discontinuous along element faces along particular interfaces, specifically material, physical, or processor boundaries.

In this paper, we present a hybrid continuous/discontinuous Galerkin spectral element method (CG/DG-SEM) for systems of linear wave equations. As a spectral element version of the CGDGTD, it has the advantages listed above for such a hybrid method, while stably allowing arbitrary order approximations on curved elements. The approximation is applicable to, and stable on, unstructured meshes on general domains with curved boundaries in up to three space dimensions with quadrilateral or hexahedral elements. We show that the approximation is stable through the use of a split form and two point fluxes. It is also conservative and constant state preserving even on meshes where the elements have curved sides.

\section{Linear Wave Propagation}

\begin{figure}[htbp] 
   \centering
   \includegraphics[width=3.5in]{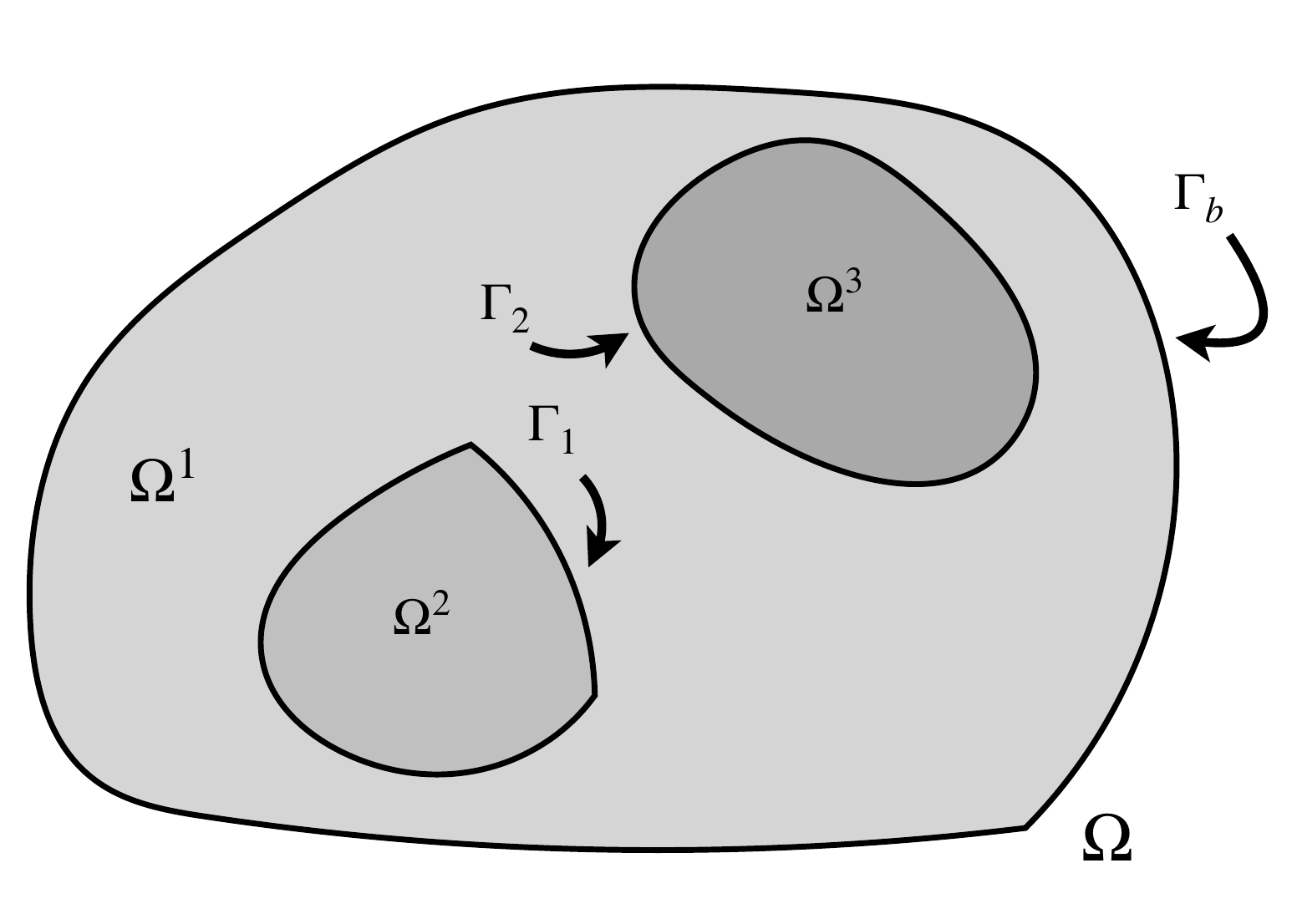} 
   \caption{Two-dimensional diagram of the domain $\Omega$ composed of $N_{mat}=3$ component subdomains, $\Omega^{m}$. Properties may jump at interior component boundaries, $\Gamma_{1}$ and $\Gamma_{2}$}
   \label{fig:MutlipleMaterialDomains}
\end{figure}

We develop the CG/DG spectral element approximation for curved hexahedral elements to solve wave propagation problems in multidimensional, multi-component domains as sketched in Fig. \ref{fig:MutlipleMaterialDomains}.
The domain $\Omega = \bigcup_{m=1}^{N_{mat}}\Omega^{m}$ is composed of ${N_{mat}}$ component subdomains, $\Omega^{m}$, each with uniform material properties and constant wave propagation speeds within each domain. Although it is possible to imbed additional components within subdomains, as $\Omega^{2}$ and $\Omega^{3}$ are embedded within $\Omega^{1}$, for simplicity of exposition we do not enable such embedding here. We are also mostly interested in the multiple material problem here, rather than parallel computing advantages, so we do not specifically consider interfaces where the inter-element properties are continuous. In any case, the numerical properties at such interfaces are already addressed in other papers, e.g. \cite{Gassner:2013ij}.

We write the wave equations on a domain $\Omega$ as hyperbolic linear systems of the form
\begin{equation}
{\statevec u_t} + \sum\limits_{i = 1}^3 {\frac{{\partial \statevec f_{i}}}{{\partial {x_i}}}}  = 0,
\label{eq:CNSE}
 \end{equation}
 where $\statevec u$ is the state vector, $\statevec f_{i} = {\mmatrix A_i}\statevec u,\;i=1,2,3$ are the covariant fluxes, $\spacevec x = (x,y,z) = \left(x_{1},x_{2},x_{3}\right)$ is the space coordinate vector, and $\mmatrix A_{i},\;i=1,2,3$ are the covariant coefficient matrices.  The coefficient matrices are expected to be piecewise constant, representing uniform properties, and are assumed to be discontinuous at component boundaries in the domain.
 
 Equations that can be written in the form \eqref{eq:CNSE} include the equations of linear acoustics \cite{PindEtAl2019},  elastic waves \cite{wilcox2010}, geometric optics \cite{VanLith2017}, and Maxwell's equations \cite{ISI:A1991GT68900017}. The example we use here is that of linear acoustics with state vector
and  coefficient matrices
\begin{equation}
\statevec u = \left[\begin{array}{c}p \\u \\v \\w\end{array}\right]
,\quad
\mmatrix A_{i} =\left[\begin{array}{cccc}0 & \delta_{i1}\rho c^2 & \delta_{i2}\rho c^2 & \delta_{i3}\rho c^2 \\\delta_{i1}/\rho & 0 & 0 & 0 \\\delta_{i2}/\rho & 0 & 0 & 0 \\\delta_{i3}/\rho & 0 & 0 & 0\end{array}\right], \quad i = 1,2,3,
\label{eq:WaveEqnMatrices}
\end{equation}
where $\rho$ is the density of the material medium, $c$ is the sound speed, and $\delta_{ij}$ is the Kronecker delta.

For the approximation to be stable, the coefficient matrices, $\mmatrix A_{i}$, must be simultaneously symmetrizable. To that end, we assume that there exists a constant symmetrizing matrix $0 < \mmatrix S < \infty$ such that $\mmatrix A_{i}^{s}=\mmatrix S^{-1} \mmatrix A_{i} \mmatrix S =  \left(\mmatrix A_{i}^{s}\right)^{T}$. For example, a symmetrization matrix for the wave equation coefficient matrices \eqref{eq:WaveEqnMatrices} is
\begin{equation}
\mmatrix S = \left[\begin{array}{cccc}c & 0 & 0 & 0 \\0 & 1/\rho & 0 & 0 \\0 & 0 & 1/\rho & 0 \\0 & 0 & 0 & 1/\rho\end{array}\right].
\end{equation}

 To compact the equations and the operations on them, we define \textit{block vectors} (with a double arrow). The flux and solution gradient, which are three-space vectors of state vectors, are block vectors,
\begin{equation}
\bigstatevec{f} =
 \left[ {\begin{array}{*{20}{c}}
  {{\statevec f_1}} \\ 
  {{\statevec f_2}} \\ 
  {{\statevec f_3}} 
\end{array}} 
\right], \quad 
\vecNablaX \statevec u = \left[ {\begin{array}{*{20}{c}}
  {{\statevec u_x}} \\ 
  {{\statevec u_y}} \\ 
  {{\statevec u_z}} 
\end{array}} 
\right].
\end{equation}
Block vectors are multiplied by block matrices. 
For example a diagonal block matrix and full block matrix are
\begin{equation}
\bigmatrix S = \left[ {\begin{array}{*{20}{c}}
  {{\mmatrix S}}\;&0&0 \\ 
  0&{{\mmatrix S}}\;&0 \\ 
  0&0&{{\mmatrix S}}\; 
\end{array}} \right]\,,\quad 
\bigmatrix B = \left[ {\begin{array}{*{20}{c}}
  {{\mmatrix B_{11}}}&{{\mmatrix B_{12}}}&{{\mmatrix B_{13}}} \\ 
  {{\mmatrix B_{21}}}&{{\mmatrix B_{22}}}&{{\mmatrix B_{23}}} \\ 
  {{\mmatrix B_{31}}}&{{\mmatrix B_{32}}}&{{\mmatrix B_{33}}} \\ 
\end{array}} \right].
\end{equation}

The vector dot product of two block vectors is a scalar and the dot product of a block vector with a space vector is a state vector, 
\begin{equation}
\bigstatevec f \cdot \bigstatevec g = \sum\limits_{i = 1}^3 {{{\statevec f}_i}^T{{\statevec g}_i}},\quad \spacevec g\cdot\bigstatevec f  = \sum\limits_{i = 1}^3 {{{ g}_i}{{\statevec f}_i}}\, .
\end{equation}
In block vector notation, then, the divergence of a flux is represented as
\begin{equation}
\vecNablaX  \cdot \bigstatevec f = \sum\limits_{i = 1}^3 {\frac{{\partial {\statevec f_i}}}{{\partial {x_i}}}}\, ,
\end{equation}
which allows the system \eqref{eq:CNSE} to be written compactly as 
\begin{equation}
  {{\statevec u}_t} + {\vecNablaX} \cdot \bigstatevec f = 0 .
  \label{eq:EqnsInBlockForm_{DAK}}
\end{equation}

We re-write  \eqref{eq:EqnsInBlockForm_{DAK}} in a split form by averaging it with the nonconservative system formed by applying the product rule to the divergence of the flux, giving
\begin{equation}
\statevec u_{t} 
+ \frac{1}{2}\left\{{\vecNablaX} \cdot \bigstatevec f+\left({\vecNablaX} \cdot \bigstatevec A\right)\statevec u 
+ \bigstatevec A\cdot{\vecNablaX}\statevec u \right\} 
= 0,
\label{eqSplitFormLinear_{DAK}}
\end{equation}
where $\bigstatevec A$ is the block vector of coefficient matrices,
\begin{equation}\bigstatevec A = \left[ {\begin{array}{*{20}{c}}\
  {{\mmatrix A_1}} \\ 
  {{\mmatrix A_2}} \\ 
  {{\mmatrix A_3}} 
\end{array}} \right].
\end{equation} 

Split forms are needed to prove stability of the approximation \cite{Gassner:2013ij}. Since the coefficient matrices can have jump discontinuities at subdomain boundaries, derivatives must of course be interpreted weakly.
Since we assume in this paper that the coefficient matrices are piecewise constant, with changes only at the component boundaries, \eqref{eqSplitFormLinear_{DAK}} reduces to
\begin{equation}
\statevec u_{t} 
+ \frac{1}{2}\left\{{\vecNablaX} \cdot \bigstatevec f+\bigstatevec A\cdot{\vecNablaX}\statevec u  \right\} 
= 0.
\label{eqSplitFormLinear_{DAK}ConstA}
\end{equation}

To completely pose the problem we must supply initial conditions, $\statevec u(\spacevec x, t) = \statevec u_{0}$ in $\Omega$, boundary conditions $\statevec g(\spacevec x,t)$ along incoming characteristics on the physical boundary, $\Gamma_{b}$, and interface conditions on the $N_{\rm intf}$ subdomain interfaces, 
\begin{equation}
\Gamma = \cup_{i}^{N_{\rm intf}}\Gamma_{i},
\end{equation}
 see Fig. \ref{fig:MutlipleMaterialDomains}. 

Conservative or nonconservative interface conditions can be specified at the component interfaces \cite{La-Cognata:2016ng},\cite{kopriva2020stability}. Conservative interface conditions require that the normal flux be continuous. Using the jump notation $\jump{u} = u_{R}-u_{L}$, where left (L) and right (R) are defined relative to a normal direction, $\spacevec n$,
the conservative interface condition is 
\begin{equation}
\jump{\bigstatevec A\cdot\spacevec n\statevec u}=0.
\label{eq:ConservJumpCondition}
\end{equation}

\begin{remark}
Sometimes it is possible to write the system \eqref{eq:CNSE} in the non-conservative form 
\begin{equation}
\mmatrix B\statevec u_{t}+ \bigstatevec {\mmatrix A}\cdot\nabla \statevec u = 0,
\label{eq:NonconservativePDE}
\end{equation}
where $\mmatrix B > 0$ and $\bigstatevec A\cdot\spacevec n$ is continuous across the subdomain interfaces. Examples include Maxwell's equations \cite{Hesthaven:2002uq}, or elastic wave equations \cite{wilcox2010}. In that case, the nonconservative interface condition is $\jump{\statevec u} = 0$.
\end{remark}

In light of the possibility that the equation form \eqref{eqSplitFormLinear_{DAK}ConstA} can be conservative or nonconservative depending on the interface condition \eqref{eq:ConservJumpCondition} or $\jump{\statevec u} = 0$, we define the conservative approximation as
\begin{definition}
The conservative scheme approximates the system \eqref{eqSplitFormLinear_{DAK}ConstA} with the conservative interface condition \eqref{eq:ConservJumpCondition}.
\label{def:ConservativeDefinition}
\end{definition}

\subsection{Energy Behavior of the Continuous Solution}

We show stability of the approximation by discretely mimicking the energy behavior of the PDE system \cite{Nordstrom:2016jk}, so it is necessary to specify the behavior to match. We define the energy through the $L_{2}$ inner product and its induced norm
\begin{equation}
\iprod{\statevec q,\statevec v}_{\Omega}=\int_{\Omega}\statevec q^{T}\statevec v dV,\quad \inorm{\statevec u}_{\Omega}=\iprod{\statevec u, \statevec u}_{\Omega}^{1/2}.
\end{equation}
 Except when it is necessary to explicitly designate the domain over which the inner product or norm is applied, we will drop the subscript $\Omega$.

To form the equation satisfied by the energy, we multiply \eqref{eqSplitFormLinear_{DAK}ConstA} by $\left(\mmatrix S^{-1}\right)^{T}\mmatrix S^{-1}\statevec u$ and integrate over $\Omega$. Then
\begin{equation}
\begin{split}
\iprod{\statevec u_{t},\left(\mmatrix S^{-1}\right)^{T}\mmatrix S^{-1}\statevec u} 
+ \oneHalf\iprod{\vecNablaX\cdot\blockvec f,\left(\mmatrix S^{-1}\right)^{T}\mmatrix S^{-1}\statevec u}
+ \oneHalf\iprod{ \bigstatevec A\cdot{\vecNablaX}\statevec u,\left(\mmatrix S^{-1}\right)^{T}\mmatrix S^{-1}\statevec u} = 0.
\end{split}
\label{eq:WellPosedStep1}
\end{equation}
Let us then define the symmetry variables as $\statevec u^{s} = \mmatrix S^{-1}\statevec u$ to write \eqref{eq:WellPosedStep1} as
\begin{equation}
\oneHalf\frac{d}{dt}\inorm{\statevec u^{s}}^{2}
+ \oneHalf\iprod{\vecNablaX\cdot\blockvec f^{s},\statevec u^{s}}
+ \oneHalf\iprod{ \bigstatevec A^{s}\cdot{\vecNablaX}\statevec u^{s},\statevec u^{s}} = 0,
\label{eq:ContEnergy0}
\end{equation}
where  $\mmatrix A_{i}^{s}= \mmatrix S^{-1}\mmatrix A_{i}\mmatrix S$, and $\statevec f_{i}^{s}=\mmatrix A_{i}^{s}\statevec u^{s}$.

To manage the discontinuities at the subdomain boundaries, we subdivide the integrals in \eqref{eq:ContEnergy0} into components. Let us re-define
\begin{equation}
\inorm{\statevec u^{s}}^{2} = \sum_{m=1}^{N_{mat}}\inorm{\statevec u^{s}}_{\Omega^{m}}^{2}.
\end{equation}
Then
\begin{equation}
\oneHalf\frac{d}{dt}\inorm{\statevec u^{s}}^{2}
+ \sum_{m=1}^{N_{mat}}\left\{\oneHalf\iprod{\vecNablaX\cdot\blockvec f^{s},\statevec u^{s}}_{\Omega^{m}}
+ \oneHalf\iprod{ \bigstatevec A^{s}\cdot{\vecNablaX}\statevec u^{s},\statevec u^{s}}_{\Omega^{m}} \right\} 
= 0.
\label{eq:ContEnergy01}
\end{equation}

The inner product of the second term in the braces of \eqref{eq:ContEnergy01} can be re-written into a volume and a surface contribution using multidimensional integration by parts. The divergence satisfies the multidimensional integration by parts formula
\begin{equation}
\iprod{\vecNabla \cdot \blockvec v,\statevec u}_{\partial\Omega^{m}} = \int_{\partial \Omega^{m}} \statevec u^{T}\blockvec v\cdot\spacevec n\dS - \iprod{  \blockvec v,\vecNabla\statevec u}_{\partial\Omega^{m}},
\label{eqw:MDIBP}
\end{equation}
for some $\statevec u,\statevec v$. So
\begin{equation}
\iprod{ \bigstatevec A^{s}\cdot{\vecNablaX}\statevec u^{s},\statevec u^{s}}_{\partial\Omega^{m}} =\iprod{ {\vecNablaX}\statevec u^{s},\bigstatevec A^{s}\statevec u^{s}}_{\partial\Omega^{m}}
=\int_{\partial \Omega^{m}} \statevec u^{s,T}\blockvec f^{s}\cdot\spacevec  n\dS - \iprod{  \vecNabla\cdot\blockvec f^{s},\statevec u^{s}}_{\partial\Omega^{m}}.
\label{eq:Continuous IntByPartsTerm}
\end{equation}
Therefore, the second term of \eqref{eq:ContEnergy1} cancels with the last term of \eqref{eq:Continuous IntByPartsTerm} leaving
\begin{equation}
\frac{d}{dt}\inorm{\statevec u^{s}}^{2} + \sum_{m=1}^{N_{mat}}\int_{\partial \Omega^{m}} \statevec u^{s,T}\blockvec f^{s}\cdot\spacevec n\dS
=0.
\label{eq:ContEnergy1}
\end{equation}

Next, we separate the surface integrals over the physical boundary, $\Gamma_{b}=\partial \Omega$, and the interior interfaces,
\begin{equation}
\frac{d}{dt}\inorm{\statevec u^{s}}^{2}  
+ \sum_{m=1}^{N_{mat}}\int_{\partial \Omega^{m}\backslash \partial \Omega} \statevec u^{s,T}\blockvec f^{s}\cdot\spacevec n\dS
+ \int_{\partial \Gamma_{b}} \statevec u^{s,T}\blockvec f^{s}\cdot\spacevec n\dS=0.
\label{eq:ContEnergySepBndry}
\end{equation}
The internal subdomain interfaces are counted twice in \eqref{eq:ContEnergySepBndry}, with normals in opposite directions. Gathering them, \eqref{eq:ContEnergySepBndry} becomes
\begin{equation}
\frac{d}{dt}\inorm{\statevec u^{s}}^{2} 
-\int_{ \Gamma} \jump{\statevec u^{s,T}\blockvec f^{s}\cdot\spacevec n}\dS
+ \int_{\Gamma_{b}} \statevec u^{s,T}\blockvec f^{s}\cdot\spacevec n\dS = 0,
\label{eq:timeDerivBound3}
\end{equation}
where a single normal vector is chosen between the two at each interface point.

To bound the physical boundary surface term in \eqref{eq:ContEnergySepBndry}, we use the fact that the system is hyperbolic so that 
\begin{equation}
\blockvec f^{s}\cdot\spacevec n =\left(\bigstatevec A^{s}\cdot\spacevec n\right) \statevec u^{s}\equiv\mmatrix A^{s}\statevec u^{s} = \mmatrix A^{s,+}\statevec u^{s} + \mmatrix A^{s,-}\statevec u^{s},
\label{eq:MatrixCharacteristicSplitting}
\end{equation}
where $\mmatrix A^{\pm} =\oneHalf\left( \mmatrix A  \pm \left|\mmatrix A\right|\right)$ are matrices with only positive or negative eigenvalues. If we apply boundary conditions $\statevec g$ to the incoming waves, then
\begin{equation}
\frac{d}{dt}\inorm{\statevec u^{s}}^{2} 
-\int_{ \Gamma} \jump{\statevec u^{s,T}\blockvec f^{s}\cdot\spacevec n}\dS
+ \int_{\Gamma_{b}} \statevec u^{s,T}\mmatrix A^{s,+}\statevec u^{s}\dS
= \int_{\Gamma_{b}} \statevec g^{s,T}\left|\mmatrix A^{s,-}\right|\statevec g^{s}\dS.
\label{eq:timeDerivBound}
\end{equation}
The third term is non-negative, so
\begin{equation}
\frac{d}{dt}\inorm{\statevec u^{s}}^{2} 
\le \int_{ \Gamma} \jump{\statevec u^{s,T}\blockvec f^{s}\cdot\spacevec n}\dS
+ \int_{\Gamma_{b}} \statevec g^{s,T}\left|\mmatrix A^{s,-}\right|\statevec g^{s}\dS.
\label{eq:timeDerivBound2}
\end{equation}

The internal interface contribution is interesting in that it is not necessarily non-negative \cite{La-Cognata:2016ng},\cite{kopriva2020stability} if the interface condition is conservative. We will call the integrand of the jump term in the integral over the interfaces
\begin{equation}
Q\left(\statevec u_{L},\statevec u_{R}\right) = \oneHalf\jump{\statevec u^{s,T}\blockvec f^{s}\cdot\spacevec n} 
=  \oneHalf\jump{\statevec u^{s,T}\mmatrix A^{s}\statevec u^{s}} 
=-\oneHalf\left\{ \statevec u_{L}^{s,T}\mmatrix A^{s}_{L}\statevec u^{s}_{L} -\statevec u_{R}^{s,T}\mmatrix A^{s}_{R}\statevec u^{s}_{R} \right\},
 \label{eq:jumpCondition1}
\end{equation}
where $\mmatrix A^{s} = \blockvec A^{s}\cdot\spacevec n$. Since the system is hyperbolic and since $\mmatrix A^{s}$ is symmetric, there exists a matrix $\mmatrix P$ such that $\mmatrix A^{s} = \mmatrix P \Lambda \mmatrix P^{T}$, and $\Lambda$ is the diagonal matrix of eigenvalues. Then we can define the characteristic variables \begin{equation} 
\statevec w = \mmatrix P^{T}\statevec u^{s} = \left[\begin{array}{c}\statevec w^{+} \\\statevec w^{-}\end{array}\right],
\end{equation}
where $\statevec w^{+}$ is associated with the positive eigenvalues of $\mmatrix A^{s}$ and $\statevec w^{-}$ is associated with the negative ones.
The characteristic variables are chosen upwind at the interface according to
\begin{equation}
\statevec w^{+}_{R}=\statevec w^{+}_{*},\quad \statevec w^{-}_{L}=\statevec w^{-}_{*},
\label{eq:UpwindedCharVars}
\end{equation}
where the `$*$' variables are computed so that the jump condition, \eqref{eq:ConservJumpCondition}, is satisfied. In terms of the characteristic variables, the jump condition is
 \begin{equation}
 \mmatrix P_{L} \Lambda_{L}\left[\begin{array}{c}\statevec w_{L}^{+} \\\statevec w_{*}^{-}\end{array}\right] =\mmatrix P_{R} \Lambda_{R}\left[\begin{array}{c}\statevec w_{*}^{+} \\\statevec w_{R}^{-}\end{array}\right],
 \label{eq:RHCondition1DSystem}
 \end{equation}
from which $\statevec w_{*}^{\pm}$  can be computed in terms of $\statevec w^{\pm}_{L/R}$ \cite{kopriva2020stability}. Let us also write
 \begin{equation}
 \Lambda = \left[\begin{array}{cc}\bar \Lambda^+ & 0 \\0 & \bar \Lambda^-\end{array}\right],\quad
  \Lambda^{+} = \left[\begin{array}{cc}\bar \Lambda^+ & 0 \\0 & 0\end{array}\right],\quad \Lambda^{-} = \left[\begin{array}{cc}0 & 0 \\0 & \bar \Lambda^-\end{array}\right].
 \end{equation} 
Then we can write $Q$ in terms of the characteristic variables as, \cite{kopriva2020stability},
\begin{equation}
Q\left(\statevec w_{L}, \statevec w_{R}\right)
 = -\oneHalf\left\{\statevec w^{+,T}_{L} \bar\Lambda^{+}_{L}\statevec w^{+}_{L} - \statevec w^{+,T}_{*} \bar\Lambda^{+}_{R}\statevec w^{+}_{*}\right\}
 - \oneHalf\left\{  \statevec w^{-,T}_{R}\left|\bar\Lambda^{-}_{R}\right|\statevec w^{-}_{R} - \statevec w^{-,T}_{*} \left|\bar\Lambda^{-}_{L}\right|\statevec w^{-}_{*}\right\}.
\label{eq:QForSystem}
\end{equation}

Finally, then, the behavior of the energy of the system with time is given by
\begin{equation}
\frac{d}{dt}\inorm{\statevec u^{s}}^{2} 
\le 2\int_{ \Gamma} Q\dS
+ \int_{\Gamma_{b}} \statevec g^{s,T}\left|\mmatrix A^{s,-}\right|\statevec g^{s}\dS.\quad\mathit{Conservative}
\label{eq:timeDerivBoundConservative}
\end{equation}

It is not necessarily true that $Q\le 0$, so it is not necessarily true that the $L_{2}$ energy norm is bounded in terms of the initial and boundary data. In fact, the energy can grow by a bounded amount from its initial value even if $\statevec g = 0$ \cite{kopriva2020stability}. However, there exists an alternate norm in which the energy is bounded \cite{La-Cognata:2016ng},\cite{kopriva2020stability}, and since the $L_{2}$ norm is easier to work with, we will follow \cite{kopriva2020stability} and use \eqref{eq:timeDerivBoundConservative} as a surrogate to determine when the CG/DG approximation is stable.

\begin{remark}
Alternatively, if the nonconservative form of the equation, \eqref{eq:NonconservativePDE}, and the non-conservative continuity condition $\jump{\statevec u^{s}}=0$ is used, then the energy is bounded in the norm $\inorm{\statevec u^{s}}_{B}^{2}= \iprod{\statevec u^{s},\mmatrix B\statevec u^{s}}$ for homogeneous boundary conditions. For along with $\jump{\mmatrix A}=0$, $Q=0$ . Following the same procedure used to find \eqref{eq:timeDerivBoundConservative},
\begin{equation}
\frac{d}{dt}\inorm{\statevec u^{s}}_{B}^{2} 
\le
 \int_{\Gamma_{b}} \statevec g^{s,T}\left|\mmatrix A^{s,-}\right|\statevec g^{s}\dS \quad\mathit{Nonconservative\; interface}.
\label{eq:timeDerivBoundNonConservative}
\end{equation}
Which, for $\statevec g = 0$, implies that $\inorm{\statevec u^{s}}_{B}\le \inorm{\statevec u_{0}^{s}}_{B}$.
\end{remark}

Finally, $\inorm{\statevec u^{s}} = \iprod{\statevec u,\left(\mmatrix S^{-1}\right)^{T}\mmatrix S^{-1}\statevec u}^{\oneHalf}$  defines an equivalent norm to the energy norm $\inorm{\statevec u}$, so we can equivalently use $\statevec u^{s}$ or $\statevec u$ when determining behavior of the energy.


\section{The CG/DG-SEM for Hyperbolic Systems}
To form the equations to be approximated, we multiply the split form equation \eqref{eqSplitFormLinear_{DAK}} by a test function $\testfuncOne\in\LTwo{\Omega}$ and integrate over the domain
\begin{equation}
\iprod{\statevec u_{t},\testfuncOne}_{\Omega} 
+ \oneHalf\iprod{\vecNablaX\cdot\blockvec f,\testfuncOne}_{\Omega}
+ \oneHalf\iprod{ \bigstatevec A\cdot{\vecNablaX}\statevec u,\testfuncOne}_{\Omega} = 0.
\label{eq:WeakForm0}
\end{equation}
For the hybrid scheme, the test functions will be further restricted to subspaces of $\LTwo{\Omega}$ depending on the continuity of the coefficient matrices. The split form is needed to ensure stability on curvilinear meshes \cite{Gassner:2013ij},\cite{Kopriva2016274}.

Normally at this point one would apply integration by parts to separate the boundary and interior contributions in \eqref{eq:WeakForm0} prior to forming a CG approximation. Alternatively, for a DG approximation one would subdivide the domain into elements and then apply integration by parts to each element's contribution. Instead, for the hybrid approximation, we first decompose the integrals over the full domain into the sum of contributions over each component subdomain, since the domain is composed of subdomains with different properties. Thus, we write
\begin{equation}
\sum_{m=1}^{N_{mat}}\left\{\iprod{\statevec u_{t},\testfuncOne}_{\Omega^{m}} + \oneHalf\iprod{\vecNablaX\cdot\blockvec f,\testfuncOne}_{\Omega^{m}}
+ \oneHalf\iprod{ \bigstatevec A\cdot{\vecNablaX}\statevec u,\testfuncOne}_{\Omega^{m}}\right\} = 0,
\label{eq:WeakForm1}
\end{equation}
which we re-write as 
\begin{equation}
\sum_{m=1}^{N_{mat}}\left\{\iprod{\statevec u_{t},\testfuncOne}_{\Omega^{m}} + \statevec r^{m}\right\}= 0.
\label{eq:WeakForm1.5}
\end{equation}
To avoid approximating functions that are discontinuous within elements by polynomials, we require that subdomains be meshed so that subdomain boundaries are element boundaries. Since only the boundaries of the subdomains will contain discontinuities, we see that $\testfuncOne \in\mathbb H^{1}\left(\Omega^{m}\right) \bigcap\LTwo{\Omega}$.

For each subdomain, we separate the boundary and interior contributions in \eqref{eq:WeakForm1} by applying  multidimensional integration by parts.  The second term becomes
\begin{equation}
\iprod{\vecNablaX\cdot\blockvec f,\testfuncOne}_{\Omega^{m}} = \int_{\partial\Omega^{m}}\testfuncOne^{T}\blockvec f\cdot\spacevec n \dS - \iprod{\blockvec f,\vecNablaX\testfuncOne}_{\Omega^{m}}.
\label{eq:FirstTermIBP}
\end{equation}
Similarly, the third term can be re-written as
\begin{equation}
\iprod{ \bigstatevec A\cdot{\vecNablaX}\statevec u,\testfuncOne}_{\Omega^{m}} = \iprod{ {\vecNablaX}\statevec u,\bigstatevec A^{(T)}\testfuncOne}_{\Omega^{m}} = \int_{\partial\Omega^{m}}\testfuncOne^{T}\blockvec f\cdot\spacevec n \dS - \iprod{\statevec u,\vecNablaX \cdot\blockvec f^{(T)}\left(\testfuncOne\right)}_{\Omega^{m}},
\label{eq:ThirdTermIBP}
\end{equation}
where
\begin{equation}
\bigstatevec A^{(T)} = \left[ {\begin{array}{*{20}{c}}\
  {{\mmatrix A^{T}_1}} \\ 
  {{\mmatrix A^{T}_2}} \\ 
  {{\mmatrix A^{T}_3}} 
\end{array}} \right]
,\quad 
\blockvec f^{(T)}\left(\testfuncOne\right) = \left[ {\begin{array}{*{20}{c}}\
  {{\mmatrix A^{T}_1}\testfuncOne} \\ 
  {{\mmatrix A^{T}_2}\testfuncOne} \\ 
  {{\mmatrix A^{T}_3}\testfuncOne} 
\end{array}} \right].
\end{equation}
When we substitute \eqref{eq:FirstTermIBP} and \eqref{eq:ThirdTermIBP} into  the components of \eqref{eq:WeakForm1} and re-order the terms,
\begin{equation}
\statevec r^{m}=
\int_{\partial\Omega^{m}}\testfuncOne^{T}\blockvec f\cdot\spacevec n  \dS
- \oneHalf \left\{ \iprod{ \statevec u,\vecNablaX\cdot\blockvec f^{(T)}(\testfuncOne) }_{\Omega^{m}} +\iprod{ \blockvec f(\statevec u),\vecNablaX \testfuncOne}_{\Omega^{m}} \right\} .
\label{eq:HyperbolicWeakFullDomain}
\end{equation}

Now we subdivide the domain into conforming hexahedral elements, which may have curved faces. Then \eqref{eq:HyperbolicWeakFullDomain} becomes the sum over all elements, which we group by subdomain in the sets $\Omega^{m}_{h},\;m=1,\ldots,N_{mat}$. There is no restriction on the number of elements per subdomain, as long as the global mesh is conforming. Thus,
\begin{equation}
\sum_{m=1}^{N_{mat}}\sum_{e\in \Omega^m_h}\left\{ \iprod{\statevec u_{t},\testfuncOne}_{e} + \statevec r^{m}_{e} \right\} = 0,
\label{eq:GlobalSum}
\end{equation}
where
\begin{equation}
\statevec r^{m}_{e} =
 \int_{\partial\Omega^{m}\cap\partial e}\testfuncOne^{T}\blockvec f\cdot\spacevec n  \dS
- \oneHalf \left\{ \iprod{ \statevec u,\vecNablaX\cdot\blockvec f^{(T)}(\testfuncOne) }_{e}
+\iprod{ \blockvec f(\statevec u),\vecNablaX \testfuncOne}_{e} \right\}
.
\end{equation}
We see, then, that an element contributes a surface flux term (like in DG) at element faces that lie on a subdomain (including physical) boundary. Away from a subdomain boundary, there are no surface terms, as in a CG approximation.

Since the test function is allowed to be independent between subdomains as in a DG scheme between elements, it follows that
\begin{equation}
\sum_{e\in \Omega^m_h}\left\{ \iprod{\statevec u_{t},\testfuncOne}_{e} + \statevec r^{m}_{e} \right\} = 0, \quad m = 1, 2, \ldots N_{mat}.
\label{eq:GlobalSumDG}
\end{equation}
Equation \eqref{eq:GlobalSumDG} implies that the approximation is stiffness summed within subdomains. Coupling between elements on subdomain boundaries is done through the surface flux as in DG. 

We evaluate the elemental contributions on the reference element $E= [-1,1]^{3}$ and so transform the inner products element by element. Elements are mapped from $E$ to an individual element $e$ by an isoparametric mapping 
\begin{equation}\label{eq:generalMapping}
\spacevec{x} = \spacevec{X_{e}}(\spacevec{\xi}),
\end{equation}
where
\begin{equation}
\spacevec{\xi} = (\xi\,,\,\eta\,,\,\zeta) = \left(\xi^1\,,\,\xi^2\,,\,\xi^3\right)
\end{equation}
are the reference space coordinates.

Now let \cite{winters2020construction}
\begin{equation}\label{eq:metric_matrix}
\bigmatrix{M} = \begin{bmatrix}
Ja_1^1 \mmatrix{I} & Ja_1^2 \mmatrix{I} & Ja_1^3 \mmatrix{I}\\[0.05cm]
Ja_2^1 \mmatrix{I} & Ja_2^2 \mmatrix{I} & Ja_2^3 \mmatrix{I}\\[0.05cm]
Ja_3^1 \mmatrix{I} & Ja_3^2 \mmatrix{I} & Ja_3^3 \mmatrix{I}\\[0.05cm]
\end{bmatrix},
\end{equation}
where \cite{Kopriva:2006er}
\begin{equation}
\label{eq:curlMetrics}
Ja_n^i = - \hat{x}_i \cdot \vecNablaXi\times \left( {{X_l}{\vecNablaXi }{X_m}} \right),\quad i = 1,2,3,\; n = 1,2,3,\;\;(n,m,l)\;\text{ cyclic},
\end{equation}
and the $J\spacevec a^{i}$ are the volume weighted contravariant basis vectors.
Then converting to reference space (c.f. \cite{winters2020construction}),
\begin{equation}\label{eq:state_divergence}
\vecNablaX  \cdot \bigstatevec{f} = \frac{1}{\mathcal J}\vecNablaXi\cdot\left({\bigmatrix{M}^T}\bigstatevec{f}\right) = \frac{1}{\mathcal J}\vecNablaXi\cdot\bigcontravec f = \frac{1}{\mathcal J}\vecNablaXi\cdot\left(\bigcontravec A \statevec u\right)
\end{equation} 
and
\begin{equation}\label{eq:state_gradient}
\vecNablaX\statevec u
=
\begin{bmatrix}
  \statevec{u}_x \\ 
  \statevec{u}_y \\ 
  \statevec{u}_z 
\end{bmatrix}
=
 \frac{1}{\mathcal J}\bigmatrix{M}
 \begin{bmatrix}
  \statevec{u}_{\xi} \\ 
  \statevec{u}_{\eta} \\ 
  \statevec{u}_{\zeta} 
\end{bmatrix}
= 
\frac{1}{\mathcal J} \bigmatrix{M}\vecNablaXi\statevec{u},
\end{equation}
where ${\mathcal J}$ is the volume Jacobian and $\bigcontravec A = \bigmatrix{M}^{T}\blockvec A$ is the block vector of the contravariant coefficient matrices.

When we represent the volume terms on the reference element, the elemental contributions become
\begin{equation}
\iprod{\statevec u_{t},\testfuncOne}_{e} = \iprod{\mathcal J\statevec u_{t},\testfuncOne}_{E}
\end{equation}
and
\begin{equation}
\statevec r^{m}_{e} = 
 \int_{\partial\Omega^{m}\cap\partial e}\testfuncOne^{T}\blockvec f\cdot\spacevec n \dS
- \oneHalf \left\{ \iprod{ \statevec u,\vecNablaXi\cdot\bigcontravec f^{(T)}(\testfuncOne) }_{E}
+\iprod{ \blockvec f(\statevec u),\bigmatrix{M}\vecNablaXi \testfuncOne}_{E} \right\}.
\end{equation}

The surface integrals can be written in either physical or reference space. The relation between physical and reference surface differentials is \cite{winters2020construction}
\begin{equation}
\dS^i = \left| J\spacevec{a}^i\right|d\xi^{j}d\xi^{k} = \mathcal J_f^i \dS_\xi^i,
\end{equation}
where the face Jacobian is $\mathcal J_{f}^i = \left|\mathcal J \spacevec{a}^i\right|$. We can write 
the normal surface flux in either reference element, $\bigcontravec f\cdot\hat{{n}}$, or physical, $\blockvec f\cdot\spacevec{n}$, variables through
\begin{equation}
\bigcontravec f\cdot\hat{{n}}^{i}\dS_\xi = \left(\boldsymbol{\bigmatrix M}^T\blockvec{f}\right)\cdot\hat{{n}}^{i}\dS_\xi = \blockvec{f}\cdot\left(\boldsymbol{\bigmatrix M}\hat{{n}}^{i}\right)\dS_\xi = \blockvec{f}\cdot\spacevec{n}\left|J\spacevec{a}^i\right|\dS_\xi = \blockvec{f}\cdot\spacevec{n}^{i}\dS.
\label{eq:dg:local-physical-fluxes-relationship}
\end{equation}
Therefore,
\begin{equation}
\int_{\partial E}\bigcontravec f\cdot\hat{n}\dS_\xi = \int_{\partial 
	e}\blockvec{f}\cdot\spacevec{n}\dS.
\label{eq:dg:surface-integrals-relation}
\end{equation}

With an abuse of notation in the surface integral, we write the elemental weak form for each element in the subdomain $\Omega^{m}$ as
\begin{equation}
\statevec r^{m}_{e} = 
 \int_{\partial\Omega^{m}\cap\partial E}\testfuncOne^{T}\bigcontravec f\cdot\hat n  \dS_{\xi}
- \oneHalf \left\{ \iprod{ \statevec u,\vecNablaXi\cdot\bigcontravec f^{(T)}(\testfuncOne) }_{E}
+\iprod{ \blockvec f(\statevec u),\bigmatrix{M}\vecNablaXi \testfuncOne}_{E} \right\}.
 \label{eq:r_{E}}
\end{equation}
The boundary contribution appears only on those element faces that lie on the subdomain boundaries. All faces internal to a subdomain are handled by the continuity of the
test functions and solution through the stiffness summation, \eqref{eq:GlobalSumDG}.

Finally, we can write \eqref{eq:r_{E}} in terms of the contravariant fluxes by moving $\bigmatrix M$ over to the other side of the second inner product of \eqref{eq:r_{E}},
\begin{equation}
\statevec r^{m}_{e} =
 \int_{\partial\Omega^{m}\cap\partial E}\testfuncOne^{T}\bigcontravec f\cdot\hat n  \dS_{\xi}
- \oneHalf \left\{ \iprod{ \statevec u,\vecNablaXi\cdot\bigcontravec f^{(T)}(\testfuncOne) }_{E}
+\iprod{ \bigcontravec f(\statevec u),\vecNablaXi \testfuncOne}_{E} \right\}.
 \label{eq:SpatialPartOnReference}
\end{equation}

\subsection{Approximation of Functions and Integrals}
We are now ready to approximate the equations, which we do by using function and integral approximations that are spectrally accurate \cite{CHQZ:2006}. 
For the spatial approximation, we represent functions of the reference domain coordinates by polynomials of degree $N$ or less, written as $\mathbb{P}^{N}(E)\subset\mathbb{L}^{2}(E)$. 

A basis for the polynomials on $E$ is the tensor product of the one dimensional Lagrange basis
\begin{equation}
{\ell _l}\left( s  \right) = \prod\limits_{i = 0;i \ne l}^N {\frac{{s  - {s _i}}}{{{s _l} - {s _i}}}}, \end{equation}
with the property $\ell_{l}\left(s_{m}\right) = \delta_{lm}$, $l,m=0,1,2,\ldots,N$. The points $s_{i},\; i=0,1,2,\ldots,N$ are the interpolation points, whose locations will be the Legendre Gauss-Lobatto quadrature points. 
 
With the Lagrange basis, we write a polynomial, $U$, in terms of its nodal values ${U}_{ijk}={U}\left(\xi_{i},\eta_{j},\zeta_{k}\right)$ as an upper case letter, which for three space dimensions is
\begin{equation}{{U}} = \sum\limits_{i,j,k = 0}^N {{{{U}}_{ijk}}{\ell _i}(\xi){\ell _j}(\eta){\ell _k}(\zeta)}. \label{eq:3DInterp}\end{equation}
One and two dimensional interpolants are special cases of \eqref{eq:3DInterp}. We also write the interpolation operator, $\mathbb{I}^{N}:\mathbb{L}^{2}(E)\rightarrow\mathbb{P}^{N}(E)$,  as
\begin{equation}{\mathbb{I}^N}\left(  u \right) = \sum\limits_{i,j,k = 0}^N {{ u\left( \xi_{i},\eta_{j},\zeta_{k}\right)}{\ell _i}(\xi){\ell _j}(\eta){\ell _k}(\zeta)}. \end{equation}

Solutions and fluxes are approximated by polynomials of degree $N$. Thus,
\begin{equation}\begin{gathered}
  {\mathbf{u}} \approx {\mathbf{U}}\in {\mathbb{P}^N}(E) \hfill \\
  \bigcontravec{\mathbf{f}} \approx \bigcontravec{\mathbf{F}}\left(\statevec U\right) = \IN{\bigcontravec A\statevec U} =\sum\limits_{i,j,k = 0}^N {{\bigcontravec A}{{\mathbf{U}}_{ijk}}{\ell _i}\left( \xi  \right){\ell _j}\left( \eta  \right){\ell _k}\left( \zeta  \right)}  \hfill \\
  \mathcal{J} \approx J = {\mathbb{I}^N}\left( \mathcal{J} \right). \hfill \\ 
\end{gathered} 
\label{eq:Approxes}
\end{equation}
Furthermore, for constant state preservation we restrict the mapping $\spacevec X_{e}$ to be isoparametric, i.e., $\spacevec X_{e}\in \PN{N}(E)$  \cite{Kopriva:2006er}. The metric terms are then approximated by polynomials of degree $N$,
\begin{equation}
\label{eq:curlMetricsP}
Ja_n^i = - \hat{x}_i \cdot \vecNablaXi\times \IN{{{X_l}{\vecNablaXi }{X_m}}},\quad i = 1,2,3,\; n = 1,2,3,\;\;(n,m,l)\;\text{cyclic}.
\end{equation}

Surface integrals and inner products are approximated by Legendre Gauss-Lobatto quadratures. 
We represent the one dimensional quadrature of a function $g(s)$ as
\begin{equation}
\int_{ - 1}^1 {g(s)ds}\approx \sum\limits_{i = 0}^N {{g\left(s_{i}\right)}{w _i}}=\sum\limits_{i = 0}^N {{g_{i}}{w _i}}\equiv \int_N {gds}    ,
 \end{equation} 
 where the $w_{i}$ are the quadrature weights,
which is exact if $g\in\mathbb{P}^{2N-1}$. By tensor product extension, three dimensional volume integral approximations are
 \begin{equation}
 \int_{E,N} {gd \xi d\eta d\zeta}  \equiv \sum\limits_{i,j,k = 0}^N {{g_{ijk}}w_{ijk}},
  \end{equation}
 where $w_{ijk}\equiv w_{i}w_{j}w_{k}$. Also by extension, we have the discrete inner product approximation
 \begin{equation}
 \iprodN{U,V} = \sum\limits_{i,j,k = 0}^N {{U_{ijk}}{V_{ijk}}w_{ijk}} \approx \iprod{U,V}.
 \end{equation}
Surface integral approximations are represented as the shorthand notation
 \begin{equation}
 \begin{split}
 \int_{\partial E,N} {\spacevec g \cdot \hat n\dS}  &= \sum\limits_{i,j = 0}^N {\left. {{w_{ij}}{g_{1}}\left( {\xi,{\eta_i},{\zeta_j}} \right)} \right|_{\xi =  - 1}^1}  + \sum\limits_{i,j = 0}^N {\left. {{w_{ij}}g_{2}\left( {{\xi_i},\eta,{\zeta_j}} \right)} \right|_{\eta =  - 1}^1}  + \sum\limits_{i,j = 0}^N {\left. {{w_{ij}}g_{3}\left( {{\xi_i},{\eta_j},\zeta} \right)} \right|_{\zeta =  - 1}^1} 
\\&\equiv\int_N {\left. {g_{1}d\eta d\zeta } \right|} _{\xi  =  - 1}^1 + \int_N {\left. {g_{2}d\xi d\zeta } \right|} _{\eta  =  - 1}^1 + \int_N {\left. {g_{3}d\xi d\eta } \right|} _{\zeta  =  - 1}^1.
 \end{split}
 \end{equation}

With the shorthand representations, the Gauss-Lobatto quadrature satisfies the important summation by parts rule \cite{10.1007/978-3-319-65870-4_2},\cite{winters2020construction}, 
 \begin{equation}
 {\iprodN{\vecNablaXi  \cdot \blockvec F,\statevec V} = \int_{\partial E ,N} {\blockvec F \cdot \hat nVdS}  - \iprodN{\blockvec F,\vecNablaXi \statevec V}},
 \label{eq:DiscreteGreens_DAK}
 \end{equation}
for any block vector of polynomials $\blockvec F\in\mathbb{P}^{N}$ and any state vector polynomial $\statevec V\in\mathbb{P}^{N}$.
Equation \eqref{eq:DiscreteGreens_DAK} is the discrete equivalent of \eqref{eqw:MDIBP}.

\subsection{The CG/DG-SEM Approximation}

We replace fluxes in \eqref{eq:SpatialPartOnReference} by polynomials and integrals by quadratures, and restrict the test functions to be continuous polynomials over each subdomain and piecewise continuous over the full domain, giving
 \begin{equation}
\statevec R^{m}_{e}=
\int_{\partial\Omega^{m}\cap\partial E,N}\testfuncOne^{T}\bigcontravec F\cdot\hat n  \dS_{\xi}
- \oneHalf \left\{ \iprodN{ \statevec U,\vecNablaXi\cdot\bigcontravec F^{(T)}(\testfuncOne) }
+\iprod{ \bigcontravec F(\statevec U),\vecNablaXi \testfuncOne}_N \right\}
.
\end{equation}
Then over each subdomain,
\begin{equation}
\sum_{e\in \Omega^m_h}\left\{ \iprodN{J\statevec U_{t},\testfuncOne} + \statevec R^{m}_{e} \right\} = 0, \quad m = 1, 2, \ldots N_{mat},
\label{eq:GlobalSumDGDiscrete}
\end{equation}
and over the whole domain,
\begin{equation}
\sum_{m=1}^{N_{mat}}\sum_{e\in \Omega^m_h}\left\{ \iprodN{J\statevec U_{t},\testfuncOne} + \statevec R^{m}_{e} \right\} = 0.
\label{eq:GlobalSumDiscrete}
\end{equation}

We also need to couple the elements across subdomain boundaries and apply boundary 
conditions to the physical boundary integrals, which is done as usual with a numerical flux. When we replace the subdomain boundary flux with the numerical flux, $\statevec F^{*}$, the CG/DG-SEM approximation of the hyperbolic system on element, $e$, is
 \begin{equation}
\begin{split}
\statevec R_{[W]} &=  \int_{\partial\Omega^{m}\cap\partial E,N}\testfuncOne^{T}\statevec F^{*} \dS_{\xi}
- \oneHalf \left\{ \iprodN{ \statevec U,\vecNablaXi\cdot\bigcontravec F^{(T)}(\testfuncOne) }
+\iprod{ \bigcontravec F(\statevec U),\vecNablaXi \testfuncOne}_N \right\},
 \end{split}
 \label{eq:RWeakForm}
\end{equation}
which is the ``weak form'' of the approximation.

The upwind numerical flux is found using the jump condition \eqref{eq:ConservJumpCondition} and the upwind evaluation of the characteristic variables, \eqref{eq:UpwindedCharVars},
\begin{equation}
\begin{split}
\statevec F^{*}\left(\statevec U_{L},\statevec U_{R}\right) &= \tilde{\mmatrix A}_{L}\mmatrix P_{L}\left[\begin{array}{c}\statevec W_{L}^{+}\\\statevec W_{*}^{-}\end{array}\right] = \tilde{\mmatrix A}_{R}\mmatrix P_{R}\left[\begin{array}{c}\statevec W_{*}^{+}\\\statevec W_{R}^{-}\end{array}\right] 
\\&
= \mmatrix P_{L}\Lambda_{L}\left[\begin{array}{c}\statevec W_{L}^{+}\\\statevec W_{*}^{-}\end{array}\right] = 
\mmatrix P_{R}\Lambda_{R}\left[\begin{array}{c}\statevec W_{*}^{+}\\\statevec W_{R}^{-}\end{array}\right],
\end{split}
\label{eq:FStarDiscontinuous}
\end{equation}
where $\tilde{\mmatrix A} = {\bigcontravec A}\cdot\hat n$ and the matrix $\mmatrix P$ now contains the eigenvectors of the matrix $ \tilde{\mmatrix A}$. 
When the coefficient matrices are constant, or at the physical boundary, the numerical flux \eqref{eq:FStarDiscontinuous} reduces to the usual form
\begin{equation}
 \statevec F^{*}\left(\statevec U_{L}, \statevec U_{R}\right)=\frac{\tilde{\mmatrix A}\statevec U_{L} + \tilde{\mmatrix A}\statevec U_{R}}{2}+\oneHalf \left| \tilde{\mmatrix A}\right|\left(\statevec U_{L} - \statevec U_{R}\right) = \tilde{\mmatrix A}\avg{\statevec U} - \oneHalf\left| \tilde{\mmatrix A}\right|\jump{\statevec U}\,,
\label{eq:LinearNumericalFlux}
\end{equation}
where  $\avg{\cdot}$ is the average of the left and right states.


Finally, like DG, it is possible to construct alternate forms using summation by parts \cite{winters2020construction}. Applying \eqref{eq:DiscreteGreens_DAK} to the two terms in braces in 
\eqref{eq:RWeakForm}, gives us an equivalent strong form
 \begin{equation}
\begin{split}
\statevec R_{[S1]} &= \int_{\partial\Omega^{m}\cap\partial E,N}\testfuncOne^{T}\statevec F^{*} \dS_{\xi}
- \int_{\partial E,N}\testfuncOne^{T}\bigcontravec F\cdot\hat n \dS_{\xi}
\\&+ \oneHalf \left\{ \iprodN{ \vecNablaXi\statevec U,\bigcontravec F^{(T)}(\testfuncOne) }
+\iprod{ \vecNablaXi\cdot\bigcontravec F(\statevec U), \testfuncOne}_N \right\}.
 \end{split}
 \label{eq:RStrongForm1}
\end{equation}
But $\statevec R_{[S1]}$ is used only in the context of the global sums, \eqref{eq:GlobalSumDGDiscrete} and \eqref{eq:GlobalSumDiscrete}. When summing over all elements, the surface fluxes from the interiors of subdomains cancel
due to continuity of the test functions and normal flux, except at the physical and subdomain boundaries. Therefore, it is acceptable to ignore those terms in the
local representation except at subdomain boundaries, leaving a {\it functionally equivalent} strong form
 \begin{equation}
\begin{split}
\statevec R_{[S]} &=  \int_{\partial\Omega^{m}\cap\partial E,N}\testfuncOne^{T}\left\{\statevec F^{*}-\bigcontravec F\cdot\hat n\right\} \dS_{\xi}
+ \oneHalf \left\{ \iprodN{ \vecNablaXi\statevec U,\bigcontravec F^{(T)}(\testfuncOne) }
+\iprod{ \vecNablaXi\cdot\bigcontravec F(\statevec U), \testfuncOne}_N \right\}.
 \end{split}
 \label{eq:RStrongForm}
\end{equation}

If, on the other hand, we only use summation by parts on the flux divergence term in \eqref{eq:RWeakForm}, and again realize that we can ignore interior surface integrals because they will cancel in the stiffness summation, then we get the ``directly-stable'' form
 \begin{equation}
\begin{split}
\statevec R_{[DS]} &=  \int_{\partial\Omega^{m}\cap\partial E,N}\testfuncOne^{T}\left\{\statevec F^{*}-\oneHalf\bigcontravec F\cdot\hat n\right\} \dS_{\xi}
+ \oneHalf \left\{ \iprodN{ \vecNablaXi\statevec U,\bigcontravec F^{(T)}(\testfuncOne) }
-\iprod{ \bigcontravec F(\statevec U), \vecNablaXi\testfuncOne}_N \right\}
.
 \end{split}
 \label{eq:RDS}
\end{equation}

\begin{remark}
We re-state for emphasis that spatial approximations $\statevec R_{[W]}$, $\statevec R_{[S]}$ and $\statevec R_{[DS]}$ are not precisely equivalent as they are for a purely DG scheme \cite{10.1007/978-3-319-65870-4_2}, but are functionally so when used in the stiffness sum \eqref{eq:GlobalSumDGDiscrete}. The surface terms left out of $\statevec R_{[S]}$ and $\statevec R_{[DS]}$ would cancel at CG interfaces, and hence contribute nothing to the global sum. Since they are not necessary, we drop those terms both in the presentation and in the implementation.
\end{remark}

\subsubsection{Two-Point Flux Form}
Each of the functionally equivalent spatial approximations, \eqref{eq:RWeakForm}, \eqref{eq:RStrongForm1}, \eqref{eq:RStrongForm}, \eqref{eq:RDS}, presented so far suggest that three derivative (matrix-vector) operations must be performed per space direction, one each to compute the volume quadratures. For example, 
\begin{equation}
\left.{\nabla \cdot\bigcontravec F} \right|_{ijk} = \sum\limits_{n = 0}^N {{\contrastatevec F^{(1)}_{njk}}{\mmatrix{D}_{in}}}  + \sum\limits_{n = 0}^N {{\contrastatevec F^{(2)}_{ink}}{\mmatrix{D}_{jn}}}  + \sum\limits_{n = 0}^N {{\contrastatevec F^{(3)}_{ijn}}{\mmatrix{D}_{kn}}}, 
\end{equation}
where $\mmatrix D_{nm} = \ell'_{n}\left( s_{m} \right)$ is the derivative matrix.
Then if we select the test function to be $\ell_{i}\ell_{j}\ell_{k}$ individually for each component (c.f. \cite{Kopriva:2009nx}),
\begin{equation}
\iprod{ \vecNablaXi\cdot\bigcontravec F(\statevec U), \ell_{i}\ell_{j}\ell_{k}}_N=  w_{ijk} \left\{\sum\limits_{n = 0}^N {{\statevec F^{(1)}_{njk}}{\mmatrix{D}_{in}}}  + \sum\limits_{n = 0}^N {{\statevec F^{(2)}_{ink}}{\mmatrix{D}_{jn}}}  + \sum\limits_{n = 0}^N {{\statevec F^{(3)}_{ijn}}{\mmatrix{D}_{kn}}}\right\}.
\end{equation}
The other volume terms can be represented similarly in such a matrix-vector multiplication form \cite{winters2020construction}.

If we add back the zero equal to the divergence of the coefficient matrices \cite{Kopriva2016274}, the three terms can be combined into a single matrix-vector differentiation of a two-point flux \cite{winters2020construction}, 
\begin{equation}
\oneHalf \left\{ \iprodN{ \vecNablaXi\statevec U,\bigcontravec F^{(T)}(\testfuncOne) }
+\iprod{ \vecNablaXi\cdot\bigcontravec F(\statevec U), \testfuncOne}_N \right\}
 + \oneHalf\iprodN{ \left({\vecNablaXi} \cdot \bigcontravec A\right)\statevec U, \testfuncOne} = \iprodN{\spacevec{\mathbb{D}} \cdot(\bigcontravec{F} )^{\#},\testfuncOne},
\end{equation}
where
the special divergence operation is
\begin{equation}
\begin{split}
\spacevec{\mathbb{D}} \cdot (\bigcontravec{F} )^{\#}(\xi,\eta,\zeta)\equiv2\sum_{n=0}^N 
&\quad \ell'_n(\xi)   \contrastatevec{F}^{\#,1}(\xi,\eta,\zeta;\xi_n,\eta,\zeta)\\[-1ex]       
&+     \ell'_n(\eta)  \contrastatevec{F}^{\#,2}(\xi,\eta,\zeta;\xi,\eta_n,\zeta)\\[1ex]    
&+     \ell'_n(\zeta) \contrastatevec{F}^{\#,3}(\xi,\eta,\zeta;\xi,\eta,\zeta_n).                                                                                                                                              
\end{split}
\label{eq:DFLinSharp}
\end{equation}
The two point flux itself can be written as the product of the average of two point values, for example,
\begin{equation}
 {{\contravec{\statevec F}}^{\#,1}_{(n,i)jk}}= {\avg{\tilde {\mmatrix A}^{1}}_{(n,i)jk}}{\avg{\statevec U}_{(n,i)jk}}.
 \label{eq:split_form_1}
\end{equation}
where
\begin{equation}
\avg{V}_{(n,i),jk}= \oneHalf\left( V_{ijk} + V_{njk} \right),
\end{equation}
etc. is the two-point average operator. Similar formulas apply to the other components.

Using the two-point flux form, we get a form equivalent to the strong form, $ \statevec R_{[S]}$ \eqref{eq:RStrongForm},
\begin{equation}
 \statevec R_{[T]} = \int_{\partial\Omega^{m}\cap\partial E,N}\testfuncOne^T\left\{\statevec F^{*}-\bigcontravec F\cdot\hat n\right\} \dS_{\xi} +\iprodN{\spacevec{\mathbb{D}} \cdot(\bigcontravec{F} )^{\#},\testfuncOne}.
\end{equation}

We gather the functionally equivalent forms of the spatial discretization inTable \ref{tab:EquivLinearForms}.
{\begin{table}[ht]
		\centering
		\caption{Functionally Equivalent CG/DG approximations to the spatial terms.}
		\label{tab:EquivLinearForms}
\begingroup
\setlength{\tabcolsep}{6pt} 
\renewcommand{\arraystretch}{2.6} 
\aboverulesep=0ex
\belowrulesep=0ex
\begin{adjustbox}{max width=\textwidth}
			\begin{tabular}{@{}l|r@{}}
				\toprule
				Form & Approximation\\
				\midrule
				\midrule
				Strong [S] & $\displaystyle\int_{\partial\Omega^{m}\cap\partial E,N}\testfuncOne^{T}\left\{\bigcontravec F^{*}-\bigcontravec F\cdot\hat n\right\} \dS_{\xi}
+ \oneHalf \left\{ \iprodN{ \vecNablaXi\statevec U,\bigcontravec F^{(T)}(\testfuncOne) }
+\iprod{ \vecNablaXi\cdot\bigcontravec F(\statevec U), \testfuncOne}_N \right\}
$ \\
				Weak [W]& $  \displaystyle\int_{\partial\Omega^{m}\cap\partial E,N}\testfuncOne^{T}\bigcontravec F^{*} \dS_{\xi}
- \oneHalf \left\{ \iprodN{ \statevec U,\vecNablaXi\cdot\bigcontravec F^{(T)}(\testfuncOne) }
+\iprod{ \bigcontravec F(\statevec U),\vecNablaXi \testfuncOne}_N \right\}
 $\\
				Directly Stable [DS] & $\displaystyle\int_{\partial\Omega^{m}\cap\partial E,N}\testfuncOne^{T}\left\{\bigcontravec F^{*}-\oneHalf\bigcontravec F\cdot\hat n\right\} \dS_{\xi}
+ \oneHalf \left\{ \iprodN{ \vecNablaXi\statevec U,\bigcontravec F^{(T)}(\testfuncOne) }
-\iprod{ \bigcontravec F(\statevec U), \vecNablaXi\testfuncOne}_N \right\}
$ 
\\
Two-Point [T] & $ \displaystyle\int_{\partial\Omega^{m}\cap\partial E,N}\testfuncOne^T\left\{\bigcontravec F^{*}-\bigcontravec F\cdot\hat n\right\} \dS_{\xi} +\iprodN{\spacevec{\mathbb{D}} \cdot(\bigcontravec{F} )^{\#},\testfuncOne} 

$
\\
				\bottomrule
		\end{tabular}
		\end{adjustbox}
		\endgroup
	\end{table}}

\section{Properties of the Approximation}
The approximation with elemental contributions listed in Table \ref{tab:EquivLinearForms} is stable. Furthermore, the conservative approximation, as defined in Def. \ref{def:ConservativeDefinition}, is globally conservative and constant state preserving for curved hexahedral elements.
\subsection{Stability}
To show stability of the CG/DG-SEM, we set $\testfuncOne = \left(\mmatrix S^{-1}\right)^{T}\mmatrix S^{-1}\statevec U$ in the form [DS], \eqref{eq:RDS}. The volume terms immediately cancel, 
because
\begin{equation}
\begin{split}
 \iprodN{\vecNablaXi\statevec U,\bigstatevec F^{\,(T)}\left( \left(\mmatrix S^{-1}\right)^{T}\mmatrix S^{-1}\statevec U  \right)} &=
 \iprodN{\bigmatrix S\vecNablaXi\mmatrix S^{-1}\statevec U,{\left(\bigmatrix S^{-1}\bigstatevec F\left( \statevec U^{s}  \right)\right)^{T}}} 
\\&= \iprodN{\vecNablaXi\mmatrix S^{-1}\statevec U,{\left(\bigmatrix S^{-1}\bigstatevec F\left( \statevec U^{s}  \right)\bigmatrix S\right)^{T}}}
\\&= \iprodN{\vecNablaXi\statevec U^{s},\bigstatevec F^{s}\left( \statevec U^{s}  \right)},
\end{split}
\label{eq:ExactFluxCancellation}
 \end{equation}
leaving only
the subdomain boundary terms. The elemental contribution to the subdomain sum \eqref{eq:GlobalSumDGDiscrete} is therefore
\begin{equation}
\iprodN{J\statevec U^{s}_{t},\statevec U^{s}} +\statevec R^{m}_{e} = \oneHalf\frac{d}{dt}\inorm{\statevec U^{s}}^{2}_{J,N}
+ \int_{\partial\Omega^{m}\cap\partial E,N}\statevec U^{s,T}\left\{\statevec F^{s,*} - \oneHalf \bigcontravec F^{s}\cdot\hat n\right\}\dS_{\xi} 
\label{eq:ElementalEnergy}
\end{equation}
where $\inorm{\statevec U^{s}}_{J,N}^{2} = \iprodN{J\statevec U^{s},\statevec U^{s}}$.

Therefore, on each subdomain,
\begin{equation}
\oneHalf\frac{d}{dt}\sum_{e\in \Omega^m_h}\inorm{\statevec U^{s}}^{2}_{J,N} +\sum_{e\in \Omega^m_h}\int_{\partial\Omega^{m}\cap\partial E,N}\statevec U^{s,T}\left\{\statevec F^{s,*} - \oneHalf \bigcontravec F^{s}\cdot\hat n\right\}\dS_{\xi} =0. 
\end{equation}
When we continue and sum over all subdomains,
\begin{equation}
\oneHalf\frac{d}{dt}\sum_{m=1}^{N_{mat}}\sum_{e\in \Omega^m_h}\inorm{\statevec U^{s}}^{2}_{J,N} 
+ \sum_{m=1}^{N_{mat}}\sum_{e\in \Omega^m_h}\int_{\partial\Omega^{m}\cap\partial E,N}\statevec U^{s,T}\left\{\statevec F^{s,*} - \oneHalf \bigcontravec F^{s}\cdot\hat n\right\}\dS_{\xi} 
=0.                                                                                                                                                                                                                                                                                                                                  
\end{equation}
As in the continuous problem, the interior subdomain interfaces combine as in a DG scheme
so that
\begin{equation}
\begin{split}
\sum_{m=1}^{N_{mat}}\sum_{e\in \Omega^m_h}\int_{\partial\Omega\cap\partial E,N}\statevec U^{s,T}\left\{\statevec F^{s,*} - \oneHalf \bigcontravec F^{s}\cdot\hat n\right\}\dS_{\xi} 
= &-\sum_{\DGfaces}\int_{N}\left\{ \jump{\statevec U^{s,T}}\statevec F^{s,*} - \oneHalf \jump{\statevec U^{s,T}\bigcontravec F^{s}\cdot\hat n}\right\}\dS_{\xi}
\\&
+\sum_{\boundaryfaces}\int_{N}\statevec U^{s,T}\left\{\statevec F^{s,*} - \oneHalf \bigcontravec F^{s}\cdot\hat n\right\}\dS_{\xi}.
\end{split}
\end{equation}
Let us now define the discrete norm over the entire domain 
\begin{equation}
\inorm{\statevec U^{s}}^{2}_{N} = \sum_{m=1}^{N_{mat}}\sum_{e\in \Omega^m_h}\inorm{\statevec U^{s}}^{2}_{J,N}.
\end{equation}
Then
\begin{equation}
\begin{split}
\oneHalf\frac{d}{dt}\inorm{\statevec U^{s}}^{2}_{N}= &+ \sum_{\DGfaces}\int_{N}\left\{ \jump{\statevec U^{s,T}}\statevec F^{s,*} - \oneHalf \jump{\statevec U^{s,T}\bigcontravec F^{s}\cdot\hat n}\right\}\dS_{\xi} 
\\&-\sum_{\boundaryfaces}\int_{N}\statevec U^{s,T}\left\{\statevec F^{s,*} - \oneHalf \bigcontravec F^{s}\cdot\hat n\right\}\dS_{\xi}.
\end{split}
\label{eq:DiscreteBound0}
\end{equation}

Comparing \eqref{eq:DiscreteBound0} with \eqref{eq:timeDerivBoundConservative}, we can state: \textit{The conservative CG/DG-SEM approximation is stable if}
\begin{equation}
-\statevec U^{s,T}\left\{\statevec F^{s,*} - \oneHalf \bigcontravec F^{s}\cdot\hat n\right\} \le \oneHalf \statevec g^{s,T}\left|\mmatrix A^{s,-}\right|\statevec g^{s}
\label{eq:DiscreteBCCondition}
\end{equation}
at physical boundary points,
and the jump condition
\begin{equation}
Q_{N}\equiv\jump{\statevec U^{s,T}}\statevec F^{s,*} - \oneHalf \jump{\statevec U^{s,T}\bigcontravec F^{s}\cdot\hat n}\le Q\left(\statevec W_{L}, \statevec W_{R}\right)
\label{eq:InterfaceCondition}
\end{equation}
holds at points along the subdomain boundaries.
For if \eqref{eq:InterfaceCondition} and \eqref{eq:DiscreteBCCondition} hold, then
\begin{equation}
\frac{d}{dt}\inorm{\statevec U^{s}}^{2}_{N}
\le  2\sum_{\DGfaces}\int_{N}QdS_{\xi} 
+\sum_{\boundaryfaces}\int_{N}\statevec g^{s,T}\left|\mmatrix A^{s,-}\right|\statevec g^{s}\dS_{\xi},
\label{eq:DiscreteEnergyBoundConservative}
\end{equation}
which is the discrete equivalent of the PDE bound \eqref{eq:timeDerivBoundConservative}.

The required bounds \eqref{eq:DiscreteBCCondition} and \eqref{eq:InterfaceCondition} have been proved previously in  \cite{Gassner:2013ij} and \cite{kopriva2020stability}. The physical boundary bound, \eqref{eq:DiscreteBCCondition}, was proved in another form in \cite{Gassner:2013ij}. We include a simpler proof in Appendix \ref{ap:Appendix}, where we show that
\begin{equation}
-\statevec U^{s,T}\left\{\bigcontravec F^{s,*}-\oneHalf\bigcontravec F^{s}\cdot\hat n\right\} 
= -\oneHalf\left(\bar{\statevec  U}^{s}_{L}
- \bar{\statevec  g}^{s} \right)^{2} + \oneHalf\bar{\statevec  g}^{s,2} 
- \oneHalf\statevec U^{s}_{L} A^{s,+}\statevec U^{s}_{L}\le \oneHalf \statevec g^{s,T}\left|\mmatrix A^{s,-}\right|\statevec g^{s},
\end{equation}
where $\bar{\statevec  U} \equiv \sqrt{\left|A^{s,-}\right|}\statevec U $, and $\bar{\statevec g}^{s}$ is defined similarly.
The bound \eqref{eq:InterfaceCondition} was proved in \cite{kopriva2020stability}, which showed that
\begin{equation}
Q_{N} = Q(\statevec W_{L}, \statevec W_{R}) -  \left(  \bar{\statevec W}^{+}_{*} - \bar{ \statevec W}_{R}^{+} \right)^{2} -  \left(  \bar{\statevec W}^{-}_{*} - \bar{ \statevec W}_{L}^{-} \right)^{2}\le Q(\statevec W_{L}, \statevec W_{R}),
\end{equation}
where $\bar {\statevec W}^{\pm} = \sqrt{\left| \bar \Lambda^{\pm}\right|}\statevec W^{\pm}$.

\begin{remark}
The approximation of the nonconservative form of the system \eqref{eq:NonconservativePDE} with the nonconservative interface leads to a stronger condition like \eqref{eq:timeDerivBoundNonConservative}. In that case, the $L_{2}$ norm is replaced by the $\mmatrix B$-weighted norm, and the interface condition with continuous coefficient matrices
\begin{equation}
\jump{\statevec U^{s,T}}\statevec F^{s,*} - \oneHalf \jump{\statevec U^{s,T}\bigcontravec F^{s}\cdot\hat n}= -\oneHalf\jump{\statevec U^{s}}^{T}\left| \tilde{\mmatrix A}^{s}\right|\jump{\statevec U^{s}}\leqslant 0
\end{equation}
was shown in \cite{Gassner2018}. Then 
\begin{equation}
\frac{d}{dt}\inorm{\statevec U^{s}}^{2}_{B,N}
\le  \sum_{\boundaryfaces}\int_{N}\statevec g^{s,T}\left|\mmatrix A^{s,-}\right|\statevec g^{s}\dS_{\xi},
\end{equation}
which is the discrete equivalent of the PDE bound \eqref{eq:timeDerivBoundNonConservative}.
\end{remark}

\subsection{Conservation}
The approximation of the conservative system (Def. \ref{def:ConservativeDefinition}) is globally conservative in the sense that the total state changes only due to flux through the physical boundary.
We show conservation with form [W]. Let $\testfuncOne = 1$. Then $\vecNablaXi \testfuncOne = 0$ and $\vecNablaXi\cdot\bigcontravec F^{(T)}(\testfuncOne) = {\vecNablaXi} \cdot \bigcontravec A = 0$. Then all three volume terms vanish. Summing over all subdomains and elements, and using the fact that the numerical flux is continuous at subdomain boundaries,
\begin{equation}
\frac{d}{dt}\sum_{e}\int_{E}J\statevec Ud\Omega = \sum_{e}\int_{\partial\Omega\cap\partial E,N}\statevec F^{*} \dS_{\xi}.
\end{equation}

\subsection{Constant State Preservation}
In wave propagation problems, especially, it is important that spurious waves are not generated by the grid and element shapes. In particular, if the solution state is constant, its time derivative should be zero, provided that the constant states satisfy the Rankine-Hugoniot (conservation) conditions at the subdomain boundaries.  

We use the strong form [S] of the spatial approximation, reproduced here,
 \begin{equation}
\begin{split}
\statevec R_{[S]} &=  \int_{\partial\Omega^{m}\cap\partial E,N}\testfuncOne^{T}\left\{\bigcontravec F^{*}-\bigcontravec F\cdot\hat n\right\} \dS_{\xi}
+ \oneHalf \left\{ \iprodN{ \vecNablaXi\statevec U,\bigcontravec F^{(T)}(\testfuncOne) }
+\iprod{ \vecNablaXi\cdot\bigcontravec F(\statevec U), \testfuncOne}_N \right\}
 \end{split}
 \label{eq:RStrongForm2}
\end{equation}
to show that a constant state is preserved.  If the state is constant within a subdomain $\Omega^{m}$, $\statevec U = \statevec C^{m}$, the gradient is zero so immediately the first volume term,
\begin{equation}
\iprodN{ \vecNablaXi\statevec C^{m},\bigcontravec F^{(T)}(\testfuncOne) } = 0.
\end{equation}
 The divergence of the contravariant flux for the constant state is $(\vecNablaXi\cdot\bigmatrix M^{T})\blockvec A^{m}\statevec C^{m}$, where we have explicitly expanded the contravariant coefficient matrix for the subdomain $\Omega^{m}$. Then
 \begin{equation}
 \iprod{ \vecNablaXi\cdot\bigcontravec F(\statevec U), \testfuncOne}_N = \iprod{ (\vecNablaXi\cdot\bigmatrix M^{T})\blockvec A^{m}\statevec C^{m}, \testfuncOne}_N .
 \label{eq:volterm2}
 \end{equation}
If the metric terms are computed isoparametrically as in \eqref{eq:curlMetricsP}, then ${\vecNablaXi} \cdot \bigmatrix M=0$  \cite{Kopriva:2006er}, and the term \eqref{eq:volterm2} also vanishes. Otherwise, this term is non-zero.

Lastly, the subdomain surface terms vanish if the Rankine-Hugoniot (conservative) conditions,
\begin{equation}
\bigcontravec F^{*}\left(\statevec U^{{m_1}},\statevec U^{{m_2}}\right) = \bigcontravec A^{{m_1}}\cdot\hat n \statevec U^{{m_1}} = -\bigcontravec A^{{m_2}}\cdot\hat n \statevec U^{{m_2}}
\end{equation}
hold for the constant states in subdomains $\Omega^{{m_1}}$ and $\Omega^{{m_2}}$ on either side of the interface. For then, in $\Omega^{m_{1}}$
\begin{equation}
\bigcontravec F^{*}-\bigcontravec F\cdot\hat n = \bigcontravec A^{{m_1}}\cdot\hat n \statevec U^{{m_1}} - \bigcontravec A^{{m_1}}\cdot\hat n \statevec U^{{m_1}} = 0
\end{equation}
and similarly in $\Omega^{m_{2}}$.

Summing over all elements and subdomains, and remembering that the form [S] has dropped the non-subdomain face contributions that cancel in the stiffness summation,
\begin{equation}
\sum_{m=1}^{N_{mat}}\sum_{e\in \Omega^m_h}\iprodN{J\statevec U_{t},\testfuncOne} = 0.
\label{eq:GlobalSumDiscretetimeDerivative}
\end{equation}
The test functions are independent at points interior to an element and at subdomain interface points, so for those points, it is immediate that
$ \left. {\statevec U_{t}}\right|_{ijk} = 0$ at any such point $ijk$. Similarly, at shared points, 
\begin{equation}
\sum_{m=1}^{N_{mat}}\sum_{e\in \Omega^m_h}\iprodN{J\statevec U_{t},\testfuncOne} = \left(\sum_{\sharedpoints}J_{ijk}w_{ijk}\testfuncOne_{ijk}\right)\left. \statevec U_{t} \right|_{shared}= 0,
\end{equation}
so the time derivative is also zero at shared CG points, and the time derivative of the solution state is zero at all points in all elements.

\section{Implementation Notes}
The implementation differs from DG because of stronger coupling between the elements due to continuity of the approximate solution and test functions. It can be split between element interior points and element face points. Interior points are decoupled, for we can choose $\testfuncOne = \ell_{i}(\xi)\ell_{j}(\eta)\ell_{k}(\zeta)$ at those points. Element face points on subdomain boundaries are treated as in a DG approximation and are coupled through the numerical flux. Face points in the interior of a subdomain are coupled through stiffness summation.

Although the four forms listed in Table \ref{tab:EquivLinearForms} are functionally equivalent, we implement the two-point flux form since it requires only one matrix vector multiplication in each direction to compute the divergence. 
Using the two point flux form, we can write down the spatial part of the elemental contribution to the global sum by replacing the test function with the local nodal basis functions, giving
\begin{equation}
\statevec R_{ijk} =
 \int_{\partial\Omega^{m}\cap\partial E,N}\ell_{i}\ell_{j}\ell_{k}\left\{\bigcontravec F^{*}-\bigcontravec F\cdot\hat n\right\} \dS_{\xi} +
 \iprodN{\spacevec{\mathbb{D}} \cdot(\bigcontravec{F} )^{\#}, \ell_{i}\ell_{j}\ell_{k}}.
\end{equation}

When we evaluate the quadratures we get the explicit representations of the pointwise values of $\statevec  R_{ijk}$ \cite{Kopriva:2009nx}.
The boundary terms are non-zero only at physical or subdomain boundary nodes so we introduce a function $\chi_{ijk}$ that is one at physical boundary and subdomain interface points and zero otherwise. Then
  \begin{equation}
\begin{split}
{\statevec R}_{ijk} &=
 \left[ {\left\{\contravec{\statevec{F}}^{*}_{Njk}-{{{{\contravec {\statevec F}}_{Njk}} \cdot \hat \xi}}\right\}{\delta_{iN}}\,w_{jk} \,\chi_{Njk}
- \left\{\contravec{\statevec{F}}^{*}_{0jk}-{{{{\contravec {\statevec F}}_{0jk}} \cdot \hat \xi}}\right\}{\delta_{i0}}\,w_{jk} \,\chi_{0jk}
+ w_{ijk}\sum\limits_{n = 0}^N {2{\contravec{\statevec F}^{\#,1}_{(n,i)jk}}{{ \dmat}_{in}}} } \right] 
 \\&+ \left[ {\left\{\contravec{\statevec{F}}^{*}_{iNk}-{{{{\contravec {\statevec F}}_{iNk}} \cdot \hat \eta}}\right\}{\delta_{jN}} \,w_{ik} \,\chi_{iNk}
- \left\{\contravec{\statevec{F}}^{*}_{i0k}-{{{{\contravec {\statevec F}}_{i0k}} \cdot \hat \eta}}\right\}{\delta_{j0}} \,w_{ik} \,\chi_{i0k}
+ w_{ijk}\sum\limits_{n = 0}^N {2{\contravec{\statevec F}^{\#,2}_{i(n,j)k}}{{ \dmat}_{jn}}} } \right] 
\\& + \left[ {\left\{\contravec{\statevec{F}}^{*}_{ijN}-{{{{\contravec {\statevec F}}_{ijN}} \cdot \hat \zeta}}\right\}{\delta_{kN}} \,w_{ij} \,\chi_{ijN}
- \left\{\contravec{\statevec{F}}^{*}_{ij0}-{{{{\contravec {\statevec F}}_{ij0}} \cdot \hat \zeta}}\right\}{\delta_{0k}}\,w_{ij} \,\chi_{ij0}
+ w_{ijk}\sum\limits_{n = 0}^N {2{\contravec{\statevec F}^{\#,3}_{ij(n,k)}}{{ \dmat}_{kn}}} } \right] .
 \end{split}
 \label{eq:PointwiseLinearWithSharpFlux}
 \end{equation}

The time derivative term is
\begin{equation}
\iprod{J\statevec U_{t},\ell_{i}\ell_{j}\ell_{k}} = J_{ijk}\dot{\statevec U}_{ijk}w_{ijk},
\end{equation}
and since element interior and DG face points are decoupled, at those points
\begin{equation}
\dot {\statevec U}_{ijk}= -\frac{1}{J_{ijk}w_{ijk}}\statevec R_{ijk}.
\label{eq:InteriorR}
\end{equation}

Element faces inside a subdomain are coupled by stiffness summation. The solution value is continuous and the test function is taken to be one at a shared point and zero everywhere else. Then all 
elements sharing that point contribute both in the time derivative term and the space approximation terms. For a particular Lagrange polynomial $\ell_{*}$ that is one at a given point $(i,j,k)=*$ and zero everywhere else, 
\begin{equation}
\sum_{e}\dot {\statevec U}^{e}_{*}J^{e}_{*}w_{*} = \sum_{e} \statevec R^{e}_{*}.
\end{equation}
Since the solution value is continuous among all the elements at the point, we can factor that out of the sum on the left so
\begin{equation}
\dot U_{*}=\frac{\sum_{e} \statevec R^{e}_{*}}{\sum_{e}J^{e}_{*}w_{*}}.
\label{eq:CGR}
\end{equation}
The sum will contain contributions only from the elements that contribute to the point.

Finally, the equations \eqref{eq:InteriorR} and \eqref{eq:CGR} are integrated in time with an explicit high order Runge-Kutta method, which has a typical CFL type restriction on the step size.

\section{Examples}
We present three examples comparing the use of CG, DG and the hybrid approach for the propagation of waves through single and multiple media. In the first example, we compute waves on a Cartesian mesh with a single material and demonstrate spectral accuracy for all three approaches. The second example is computed on domain with curved elements. The final example demonstrates spectral accuracy for plane wave reflection at a material interface.

\subsection{Planewave Propagation on a Cartesian Mesh}
In the first example we solve the acoustic wave equation in the split form \eqref{eqSplitFormLinear_{DAK}} with the coefficient matrices \eqref{eq:WaveEqnMatrices} in two space dimensions on the domain $\Omega = [-5,5]^{2}$. The mesh consists of 20 elements of equal size in each direction. Time steps in all examples were chosen small enough so that the error is dominated by the spatial approximation.

The exact plane wave solution to the system is
\begin{equation}
\statevec u = a\psi\left(\spacevec k\cdot\spacevec x - \omega(t - t_{0})\right)\left[\begin{array}{c}1 \\\frac{k_x}{\rho c} \\\frac{k_y}{\rho c}\end{array}\right],
\label{eq:WavePacketSolution}
\end{equation}
where $\psi$ is a given wave function, $a$ is the amplitude, $\spacevec k$ is the wavevector, and $\omega$ is the frequency. We present results for
 $\rho = 1$ and $c = 1$, $a=1$, $\spacevec k = \left(\frac{\sqrt{3}}{2}, \oneHalf\right)$, $\omega = 5\pi/2$ and, $t_{0}=3$. For the wave function, we choose the wavepacket
 \begin{equation}
 \psi(s) = \sin(\omega s)e^{-s^{2}/(\omega\sigma)^{2}},
 \end{equation}
 with $\sigma^{2} = -(MT)^{2}/(4\ln(10^{-4}))$, where $M = 4$ is the number of modes in the significant part of the envelope and $T = 2\pi/\omega$ is the period.

\begin{figure}[htbp] 
   \centering
   \includegraphics[width=5in]{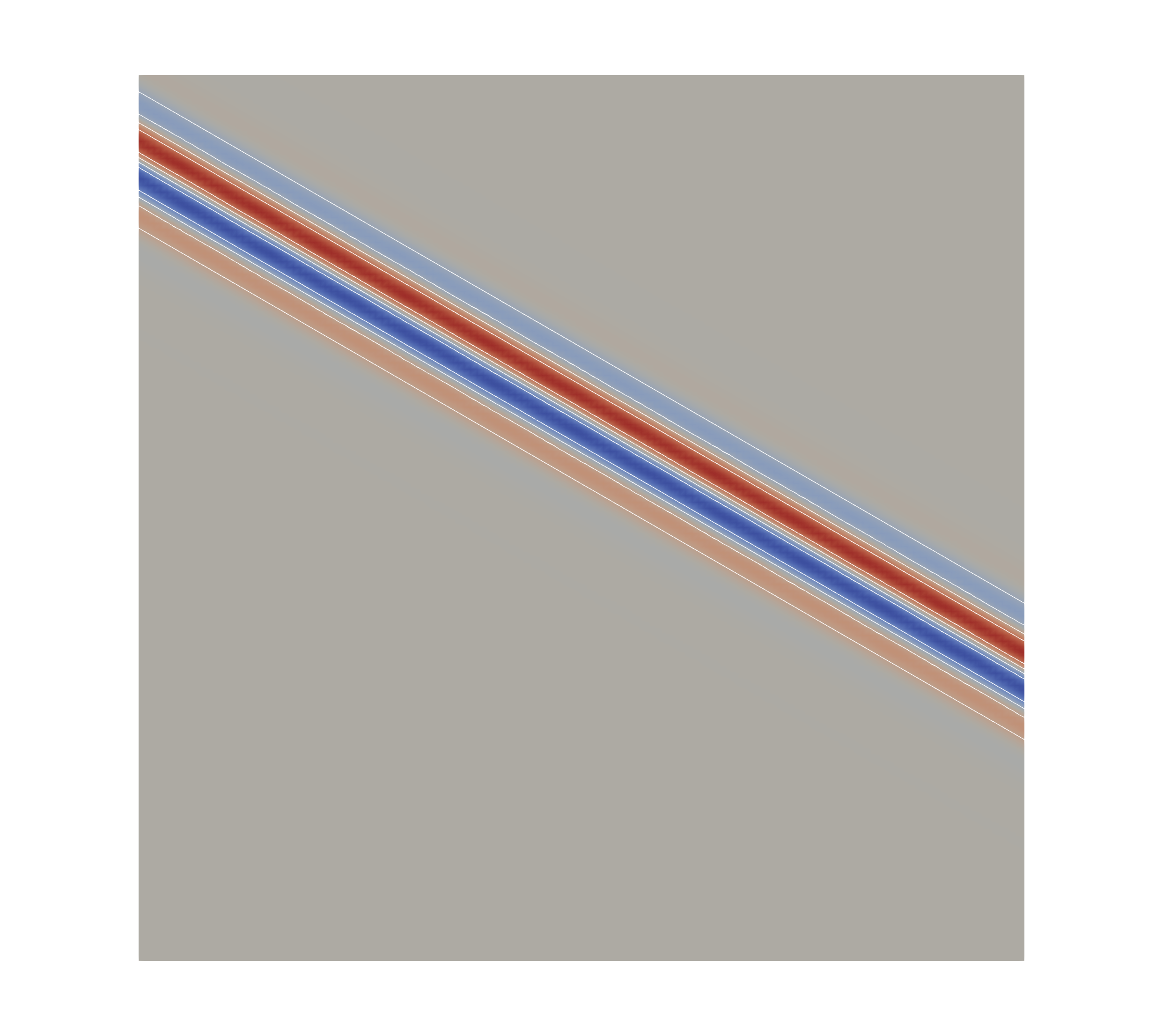} 
   \caption{Contours of $p$ for the CG approximation to \eqref{eq:WavePacketSolution} on the Cartesian mesh at time $t=4$. } 
   \label{fig:CartesianWave}
\end{figure}

\begin{figure}
    \centering
    \begin{tikzpicture}
        \begin{axis}[xlabel=$N$,ylabel=$Log_{10}(Max\; Error)$]
	\addplot[black,mark=o] coordinates {
		(4,-0.5128)
		(5,-1.6151)
		(6,-1.90655)
		(7,-3.0542)
		(8,-3.422482)
		(9,-4.668811)
            };
            \addplot[black, mark=*] coordinates {
		(4,-0.95432)
		(5,-1.523428)
		(6,-2.29684)
		(7,-3.11473)
		(8,-3.8840)
		(9,-4.69457)
	};
             \addplot[black,mark=+] coordinates {
		(4,-0.512074)
		(5,-1.61334)
		(6,-1.906541)
		(7,-3.054099)
		(8,-3.422526)
		(9,-4.66911)
	};
           \legend{CG,DG,CG/DG}
        \end{axis}
    \end{tikzpicture}
    \caption{Convergence of DG, CG and hybrid approximations for wavepacket convergence on a Cartesian mesh}
    \label{fig:CartesianErrorConvergence}
\end{figure}
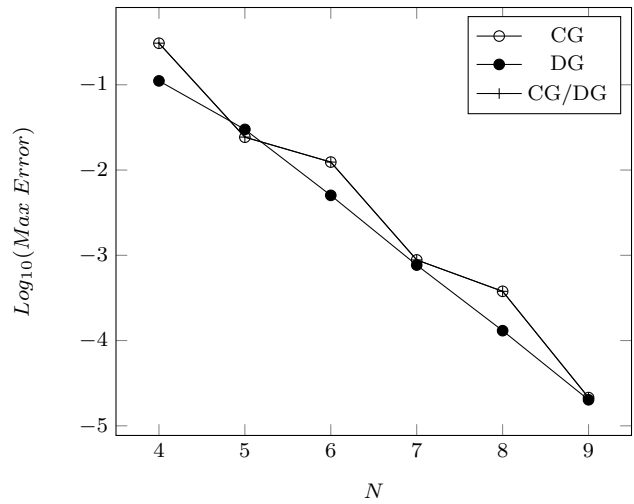

Results for pure CG, pure DG, and a CG/DG approximation are shown in Figs. \ref{fig:CartesianWave} and \ref{fig:CartesianErrorConvergence}. In the case of the hybrid approximation, all interfaces are CG except along the line $x=0$, where DG is applied, even though the properties are continuous there. Fig.  \ref{fig:CartesianWave} shows solution contours for the CG solution.

Error convergence at time $t=5$ is shown in Fig. \ref{fig:CartesianErrorConvergence}. The hybrid CG/DG mesh produces the same error as the CG to four digits. Most important to note is that the CG approximation shows an even-odd convergence behavior, whereas the full DG approximation, with an upwind numerical flux at all element interfaces, does not in this example. Note, however, that DG can also exhibit such behavior when using the central flux \cite{hindenlang2019order}. For all the approximations the error decay is exponential, and doubling the number of nodes in each direction decreases the error by about three orders of magnitude.

 \subsection{Planewave Propagation on a Curved Element Mesh}
\begin{figure}[htbp] 
   \centering
   \includegraphics[width=4in]{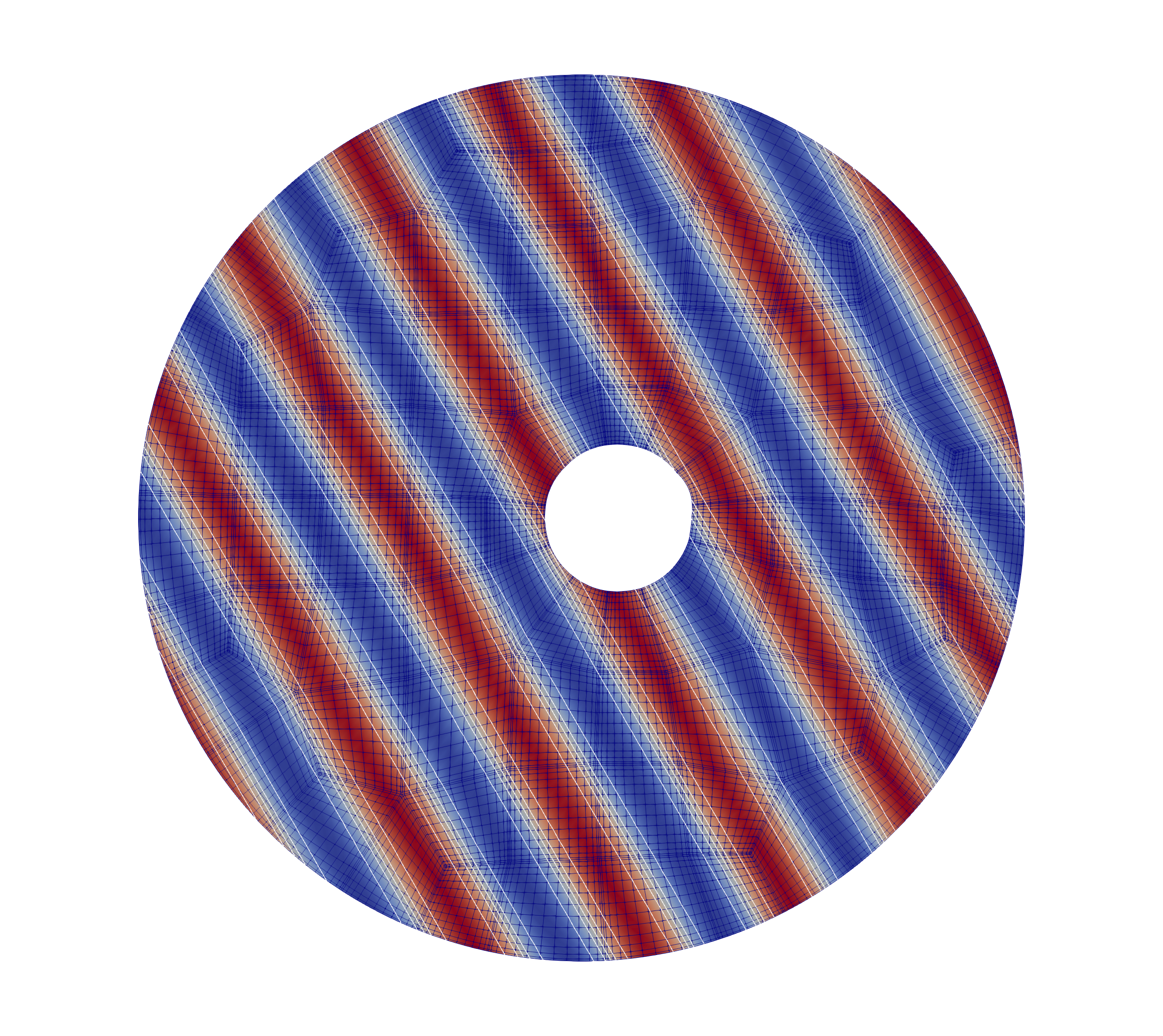} 
   \caption{Contours of $p$ for the CG approximation to a sinusoid on a non-Cartesian mesh with curved elements}
   \label{fig:CGSEMExample}
\end{figure}

 We compute the second example on a non-Cartesian mesh with curved elements. This time, we choose $\psi$ to be a sine function. Contours of $p$ from the CG approximation are shown in Fig.  \ref{fig:CGSEMExample}, overset on the spectral element mesh with internal degrees of freedom. 

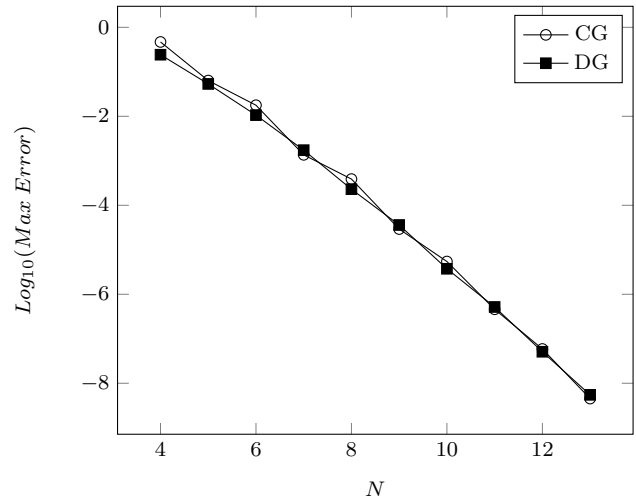
\begin{figure}
    \centering
    \begin{tikzpicture}
        \begin{axis}[xlabel=$N$,ylabel=$Log_{10}(Max\; Error)$]
	\addplot[black,mark=o] coordinates {
		(4,-0.329076)
		(5,-1.19659)
		(6,-1.75099)
		(7,-2.865967)
		(8,-3.41088769)
		(9,-4.53424)
		(10,-5.263501)
		(11,-6.337509)
		(12,-7.229493)
		(13,-8.344295)
            };
            \addplot[black,mark=square*] coordinates {
		(4,-0.61582)
		(5,-1.27228)
		(6,-1.97317)
		(7,-2.76114)
		(8,-3.6381)
		(9,-4.443067)
		(10,-5.42679)
		(11,-6.28582)
		(12,-7.2918050)
		(13,-8.259943)
	};
            \legend{CG,DG}
        \end{axis}
    \end{tikzpicture}
    \caption{Convergence of DG and CG for the non-Cartesian mesh}
    \label{fig:NoncartesianConvergence}
\end{figure}

We show the spectral error convergence on the non-Cartesian mesh in Fig. \ref{fig:NoncartesianConvergence}. As in the previous example, the CG approximation shows even-odd behavior whereas the DG does not, though not to as pronounced an effect. Doubling the polynomial order from six to 12 decreases the error by over five orders of magnitude.

\subsection{Constant State Preservation}
We demonstrate the expected constant state preservation using the curved element grid shown in Fig. \ref{fig:CGSEMExample}, where the curved circular boundaries are approximated by polynomials of order five. For that geometry, the approximation is isoparametric for $N=5$, and so the metric identities hold discretely for $N\ge 5$ \cite{Kopriva:2006er}.

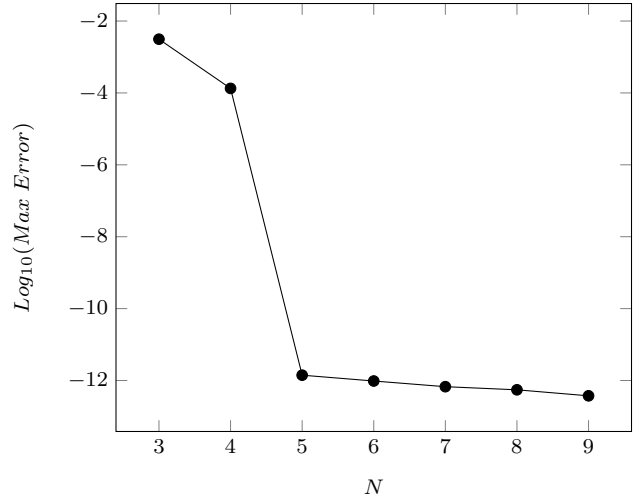
\begin{figure}
    \centering
    \begin{tikzpicture}
        \begin{axis}[xlabel=$N$,ylabel=$Log_{10}(Max\; Error)$]
	\addplot[black,mark=*] coordinates {
		(3,-2.505)
		(4,-3.874)
		(5,-11.852)
		(6,-12.014 )
		(7,-12.172)
		(8,-12.259)
		(9,-12.426)
            };
        \end{axis}
    \end{tikzpicture}
    \caption{Convergence of the CG approximation on a constant state solution and the mesh of Fig. \ref{fig:CGSEMExample}}
    \label{fig:FSPErrorConvergence}
\end{figure}

We plot the maximum error for a constant state as a function of polynomial order in Fig. \ref{fig:FSPErrorConvergence}. As expected, the error is non-zero when the boundary approximation is superparametric. As soon as the geometry and solution approximation orders match, the error drops to rounding error levels.

\subsection{Wave Scattering at a Material Interface}
\begin{figure}[htbp] 
   \centering
   \includegraphics[width=3.5in]{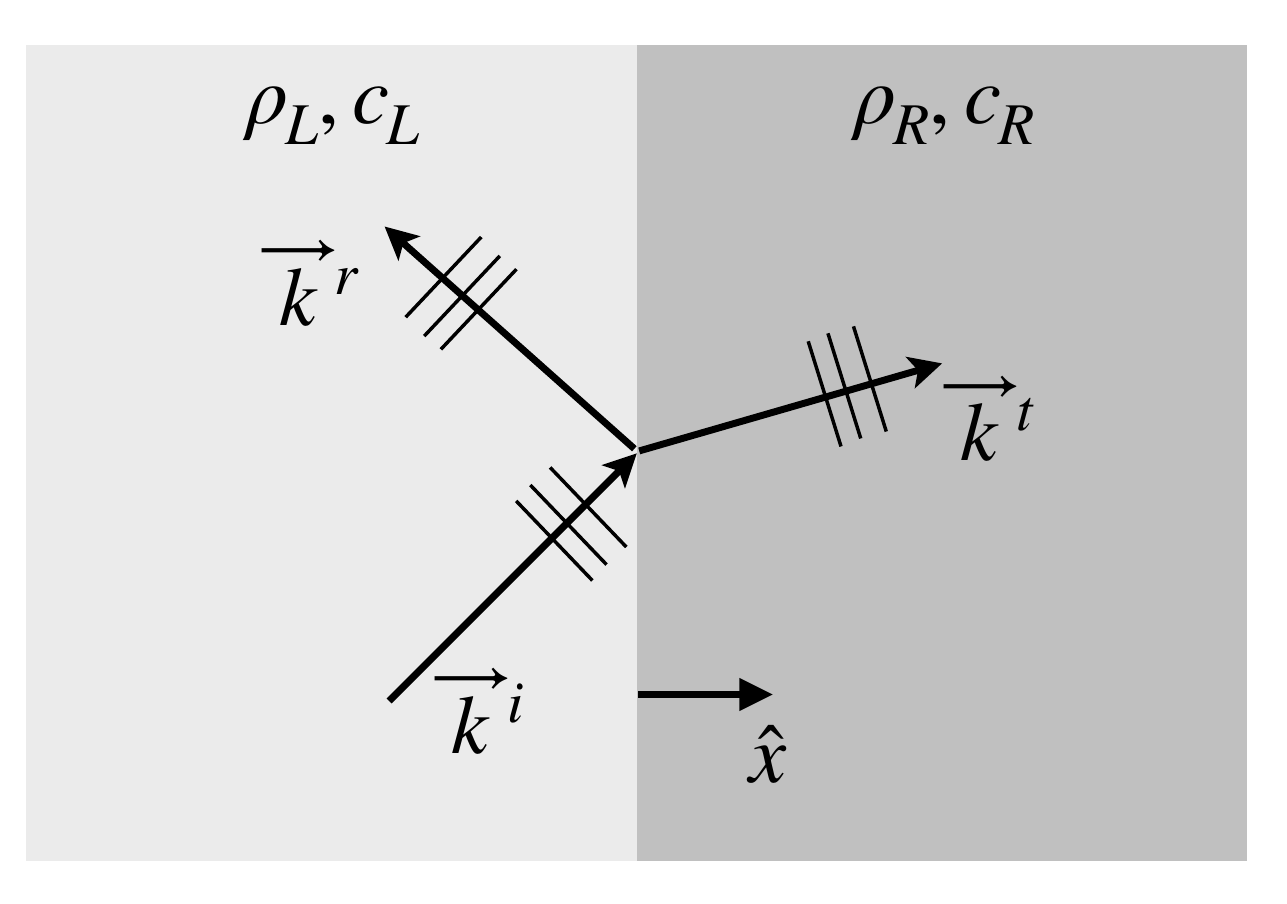} 
   \caption{Geometry of wave scattering at a plane material interface}
   \label{fig:WaveScatteringDiagram}
\end{figure}

In the final example, we show exponential convergence for the hybrid CD/DG approximation of wave scattering at a material interface as sketched in Fig. \ref{fig:WaveScatteringDiagram}. 

The exact solution is the sum of incident, reflected and transmitted waves. With the incident wavevector defined to be
\begin{equation}
\spacevec k^i  = \frac{\omega }{{c_L }}\left( {k_x^i \hat x + k_y^i \hat y} \right),
\end{equation}
the reflected and transmitted wavevectors are
\begin{equation}
\begin{array}{l}
 \spacevec k^r  = \dfrac{\omega }{{c_L }}\left( { - k_x^i \hat x + k_y^i \hat y} \right) \\ \\
 \spacevec k^t  = \dfrac{\omega }{{c_R }}\left[ {\sqrt {1 - \left( {\dfrac{{c_R }}{{c_L }}} \right)^2 \left( {k_y^i } \right)^2 } \hat x + \dfrac{{c_R }}{{c_L }}k_y^i \hat y} \right]. \\ 
 \end{array}
\end{equation}
The corresponding exact amplitudes are
\begin{equation}
\begin{array}{l}
 \dfrac{{a^r }}{{a^i }} = \dfrac{1}{J}\left( {\rho _R c_R k_x^t /k^t  - \rho _L c_L k_x^i /k^i } \right), \\ \\
 \dfrac{{a^t }}{{a^i }} = \dfrac{1}{J}\left( {\rho _L c_L k_x^r /k^r  - \rho _R c_L k_x^i /k^i } \right), \\ 
 \end{array}
\end{equation}
where $k = |k|$, and 
\begin{equation}
J =  - \rho _R c_R k_x^t /k^t  + \rho _L c_L k_x^r /k^r .
\end{equation}
Then each component of the solution has the form \eqref{eq:WavePacketSolution} with the appropriate amplitude and wavevector.

We compute the solution on the square domain $\Omega = [-5,5]^{2}$ with the material interface at $x = 0$. The mesh is Cartesian with 20 equal size elements in each direction. Other parameters for the problem are shown in Table \ref{tab:ReflTransTable}.
	 \begin{table}[tbp]
	\begin{center}
	\caption {Parameters for Plane Wave Reflection Problem}
	\label{tab:ReflTransTable}
	\begin{tabular}{c|ccccccccc}
	Parameter & $M$ & $\omega$ & $k^i_x$ & $k^i_y$ & $\rho_L$ & $\rho_R$ & $c_L$ & $c_R$ & $t_{0}$ \\
	\hline
	Value & 4 & $5\pi/2$ & 0.5 & $\sqrt{3/2}$ & 1 & 0.4 & 1 & 0.7 & 3
	\end{tabular}
	\end{center}
	\end{table}%

\begin{figure}[htbp] 
   \centering
   \includegraphics[width=4.5in]{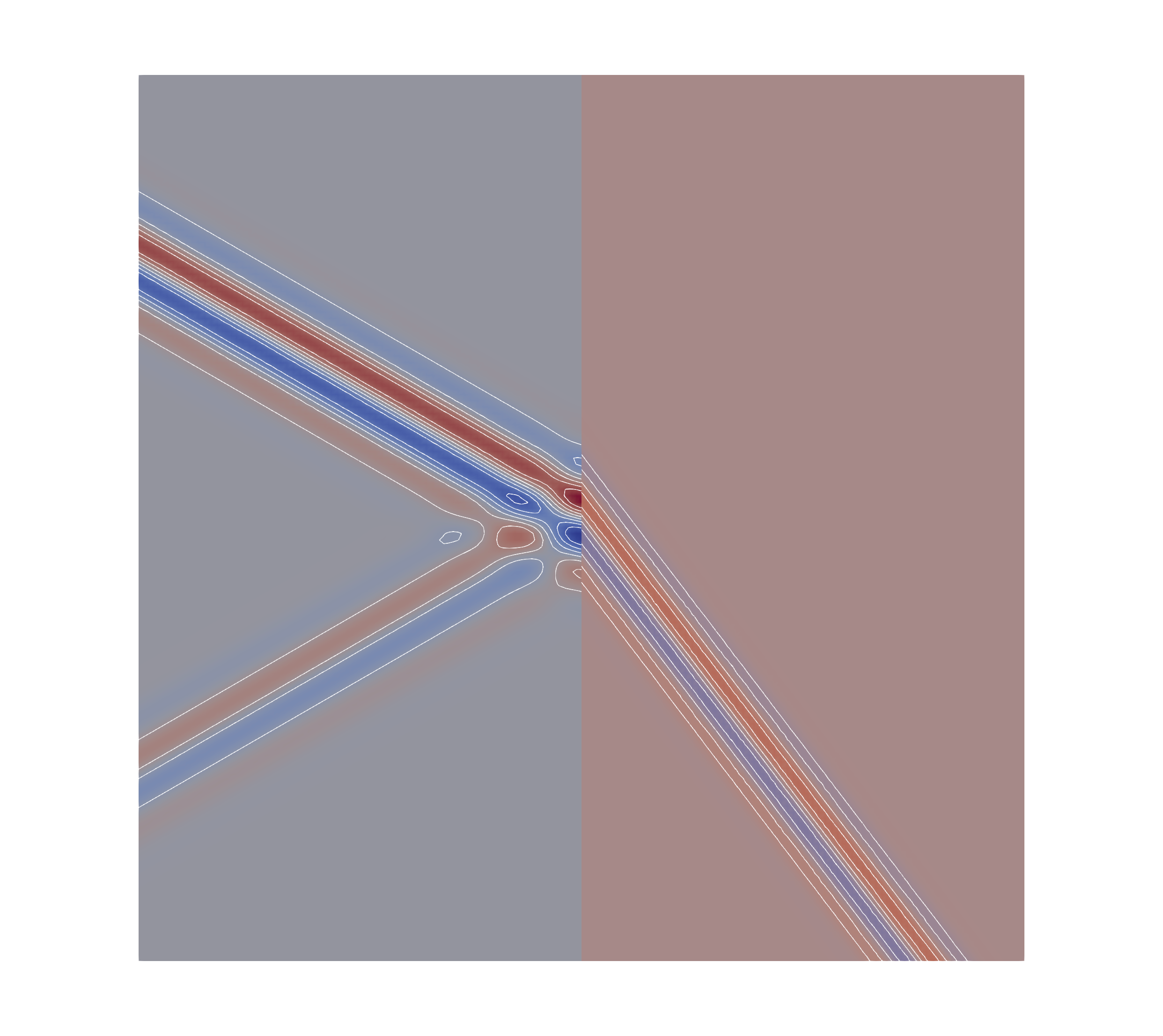} 
   \caption{Contours of $p$ for the transmission and reflection of a wavepacket at a material interface using the hybrid CG/DG approximation}
   \label{fig:ScatteringPlot}
\end{figure}

\begin{figure}
    \centering
    \begin{tikzpicture}
        \begin{axis}[xlabel=$N$,ylabel=$Log_{10}(Max\; Error)$]
	\addplot[black,mark=*] coordinates {
		(4,0.001643)
		(5,-0.77800961247999956)
		(6,-1.216607)
		(7,-1.93779)
		(8,-2.606075)
		(9,-3.19216)
		(10,-4.1106)
		(11,-4.633)
		(12,-5.45)
            };
        \end{axis}
    \end{tikzpicture}
    \caption{Convergence of the hybrid CG/DG approximation for scattering at a material interface}
    \label{fig:ScatteringErrorConvergence}
\end{figure}
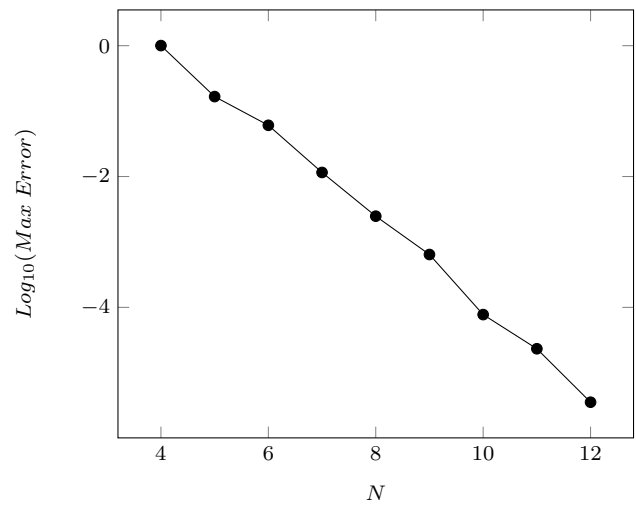

Fig. \ref{fig:ScatteringPlot} shows the contours of $p$ at time $t = 5$ for $N = 10$. Clearly seen at the material interface is the solution jump, which is where the DG approximation is applied. Exponential convergence of the error is shown in Fig. \ref{fig:ScatteringErrorConvergence}, where increasing the polynomial order from five to ten decreases the error by over three orders of magnitude. 


\begin{figure}
    \centering
    \begin{tikzpicture}
        \begin{axis}[xlabel=$t$,ylabel=$\inorm{\statevec u}_{L_{2}}$]
	\addplot+[black,only marks, mark=o] coordinates {
(0.50000000000000000,	2.9574016157116785)
(1.0000000000000000,	3.1252753792760988)
(1.5000000000000000,	3.2714852556315313)
(2.0000000000000000,	3.4040145004693434)
(2.5000000000000000,	3.5315182186018057)
(3.0000000000000000,	3.6545755782248031)
(3.5000000000000000,	3.7736216863751268)
(4.0000000000000000,	3.8888005343827170)
(4.5000000000000000,	3.9639596164674269)
(5.0000000000000000,	3.9237453162838660)
(5.5000000000000009,	3.8035099728520345)
(6.0000000000000009,	3.6759953744713152)
(6.5000000000000009,	3.5438965806356904)
(7.0000000000000009,	3.4059127321601137)
(7.5000000000000009,	3.2307084085422351)
(8.0000000000000000,	3.0623990708348376)
(8.5000000000000000,	2.8900290429694904)
(9.0000000000000000,	2.7068442715172898)
(9.5000000000000000,	2.5104650687135752)
(10.000000000000000,	2.3083686862975630)
(10.500000000000000,	2.1072722890446163)
(11.000000000000000,	1.8858162414013053)
(11.500000000000000,	1.6346289542530470)
(12.000000000000000,	1.3370637826362421)
(12.500000000000000,	0.95055016077879262)
(13.000000000000000,	0.35526218613087002)
(13.500000000000000,	1.2609492910614127E-002)
(14.000000000000000,	5.5402052575658448E-003)
(14.500000000000000,	5.2755309981154380E-003)
(15.000000000000000,	4.9859524192529649E-003)
(15.500000000000000,	4.8396322351756116E-003)
(16.000000000000000,	4.6036786332554350E-003)
            };
	\addplot+[black,no markers] coordinates {
(0.50000000000000000,	2.9574024798797254)
(1.0000000000000000,	3.1252787366683892)
(1.5000000000000000,	3.2714915990635594)
(2.0000000000000000,	3.4040241714957444)
(2.5000000000000000,	3.5315316105563141)
(3.0000000000000000,	3.6545930555480508)
(3.5000000000000000,	3.7736434128746668)
(4.0000000000000000,	3.8888268595590820)
(4.5000000000000000,	3.9639929502997449)
(5.0000000000000000,	3.9237839772476315)
(5.5000000000000009,	3.8035548614045727)
(6.0000000000000009,	3.6760454298504319)
(6.5000000000000009,	3.5439510010950470)
(7.0000000000000009,	3.4059593595133459)
(7.5000000000000009,	3.2301322608572409)
(8.0000000000000000,	3.0621132483048381)
(8.5000000000000000,	2.8899415996966886)
(9.0000000000000000,	2.7068419013740810)
(9.5000000000000000,	2.5105003957510905)
(10.000000000000000,	2.3084189091696561)
(10.500000000000000,	2.1073279128658267)
(11.000000000000000,	1.8858675968359839)
(11.500000000000000,	1.6346749429738079)
(12.000000000000000,	1.3370965371964394)
(12.500000000000000,	0.95056302362827505)
(13.000000000000000,	0.35507645608547228)
(13.500000000000000,	1.1141847886655383E-002)
(14.000000000000000,	6.6861339789815524E-006)
(14.500000000000000,	5.1531595239721890E-011)
(15.000000000000000,	0.0)
(15.500000000000000,	0.0)
(16.000000000000000,	0.0)
            };
            \legend{CG/DG,Exact}
        \end{axis}
    \end{tikzpicture}
    \caption{Time variation of the $L_{2}$ energy for the scattering problem shown in Fig. \ref{fig:ScatteringPlot}}
    \label{fig:EnertVsTime}
\end{figure}
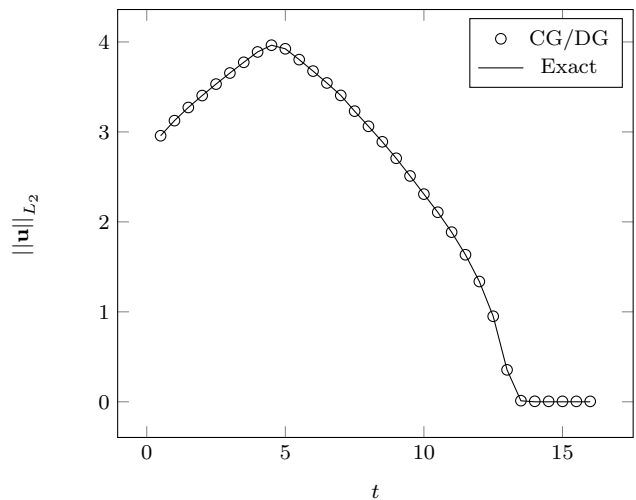

Regarding stability, we plot the energy behavior measured as the $L_{2}$ norm as a function of time in Fig. \ref{fig:EnertVsTime}. As expected from \eqref{eq:timeDerivBoundConservative} and \eqref{eq:DiscreteEnergyBoundConservative}, the energy of the continuous and discrete problem have the same behavior and is not bounded by the intitial value. Nevertheless the $L_{2}$ energy is bounded for both and vanishes as the waves exit the domain.

\section{Summary}
We have presented a stable, conservative, constant-state preserving spectral element method for isoparametric curved quadrilateral or hexahedral elements. It is a hybrid continuous and discontinuous Galerkin method where faces along subdomain boundaries use a DG numerical flux and otherwise use a continuous Galerkin ansatz. The hybrid approach allows the strengths of each method to be exploited. Numerical experiments show that the method is spectrally convergent for continuous and for discontinuous solutions when the discontinuities fall along element boundaries. 
\acknowledgement
{
 This work was supported by a grant from the Simons Foundation (\#426393, David Kopriva). Gregor Gassner thanks the Klaus-Tschira Stiftung and
the European Research Council for funding through the ERC Starting Grant “An
Exascale aware and Un-crashable Space-Time-Adaptive Discontinuous Spectral
Element Solver for Non-Linear Conservation Laws” (EXTREME, project no. 71448). }

\appendix
\section{Proof of Physical Boundary Dissipation}\label{ap:Appendix}
We show that \eqref{eq:DiscreteBCCondition} holds. As for the continuous problem, 
let $\mmatrix A^{s} = \bigcontravec A^{s}\cdot\hat n$ be the normal coefficient matrix. The physical boundary can be viewed as being between two identical media so the numerical flux reduces to the standard upwind flux \eqref{eq:LinearNumericalFlux},
\begin{equation}
\statevec F^{s,*}\left( \statevec U_{L}, \statevec U_{R}\right) =  \mmatrix A^{s,+}\statevec U_{L} + \mmatrix A^{s,-}\statevec U_{R} = \oneHalf\left( \statevec F^{s}_{L} + \statevec F^{s}_{R} \right) - \oneHalf\left| \mmatrix A^{s}\right|\left(\statevec U_{R} - \statevec U_{L}\right).
\end{equation}
Then
\begin{equation}
\statevec F^{s,*}-\oneHalf\bigcontravec F^{s}\cdot\hat n = \oneHalf \mmatrix A^{s} \left( \statevec U_{L} + \statevec U_{R}\right) - \oneHalf \left| \mmatrix A^{s}\right| \left( \statevec U_{R} - \statevec U_{L}\right) - \oneHalf \mmatrix A^{s}\statevec U_{L} = \oneHalf\left( \mmatrix A ^{s}- \left| \mmatrix A^{s}\right| \right)\statevec U_{R}+ \oneHalf\left| \mmatrix A^{s}\right|\statevec U_{L}.
\end{equation}
Since $\statevec U^{s}$ is the value taken from the interior, i.e. $\statevec U_{L}$, 
\begin{equation}
\statevec U^{s,T}\left\{\statevec F^{s,*}-\oneHalf\bigcontravec F^{s}\cdot\hat n\right\} =  \statevec U_{L}\mmatrix A^{s,-}\statevec U_{R} + \oneHalf \statevec U_{L}\left| \mmatrix A^{s}\right|\statevec U_{L}
\end{equation}
But $\left| \mmatrix A^{s}\right| = \mmatrix A^{s,+}- \mmatrix A^{s,-}$, so
\begin{equation}
\statevec U^{s,T}\left\{\statevec F^{s,*}-\oneHalf\bigcontravec F^{s}\cdot\hat n\right\} =  \statevec U_{L}\mmatrix A^{s,-}\statevec U_{R} + \oneHalf \statevec U_{L} \mmatrix A^{s,+}\statevec U_{L}-\oneHalf\statevec U_{L} \mmatrix A^{s,-}\statevec U_{L}.
\end{equation}
Also, $\mmatrix A^{s,-} = -\left| \mmatrix A^{s,-}\right|$, so
\begin{equation}
\statevec U^{s,T}\left\{\statevec F^{s,*}-\oneHalf\bigcontravec F^{s}\cdot\hat n\right\} =  -\statevec U_{L}\left|\mmatrix A^{s,-}\right|\statevec U_{R}  +\oneHalf\statevec U_{L} \left|\mmatrix A^{s,-}\right|\statevec U_{L} + \oneHalf\statevec U_{L} \mmatrix A^{s,+}\statevec U_{L} 
\end{equation}
The last two terms are non-negative.

To simplify the next few steps, we define $\bar{\statevec  U} \equiv \sqrt{\left|\mmatrix A^{s,-}\right|}\statevec U $. Then
\begin{equation}
\statevec U^{s,T}\left\{\statevec F^{s,*}-\oneHalf\bigcontravec F^{s}\cdot\hat n\right\} =    \oneHalf\left(\bar{\statevec  U}_{L}\right)^{2} -\bar{\statevec  U}_{L}\bar{\statevec  U}_{R}+ \oneHalf\statevec U_{L} \mmatrix A^{s,+}\statevec U_{L} .
\end{equation}
Completing the square,
\begin{equation}
\statevec U^{s.T}\left\{\bigcontravec F^{s,*}-\oneHalf\bigcontravec F^{s}\cdot\hat n\right\} = \oneHalf\left(\bar{\statevec  U}_{L} - \bar{\statevec  U}_{R} \right)^{2} - \oneHalf\left(\bar{\statevec  U}_{R}\right)^{2} + \oneHalf\statevec U_{L} \mmatrix A^{+}\statevec U_{L}.
\end{equation}
Replacing the external state $\bar{\statevec U}_{R}$ with a boundary condition $\bar{\statevec g}$, 
\begin{equation}
\statevec U^{s,T}\left\{\bigcontravec F^{s,*}-\oneHalf\bigcontravec F^{s}\cdot\hat n\right\} 
= \oneHalf\left(\bar{\statevec  U}_{L} 
- \bar{\statevec  g} \right)^{2} - \oneHalf\bar{\statevec  g}^{2} 
+ \oneHalf\statevec U_{L} \mmatrix A^{s,+}\statevec U_{L}.
\end{equation}
Therefore,
\begin{equation}
-\statevec U^{s,T}\left\{\bigcontravec F^{s,*}-\oneHalf\bigcontravec F^{s}\cdot\hat n\right\}  \le   \oneHalf\bar{\statevec  g}^{2}.
\end{equation}
When we return to the original variables, we get the desired result, \eqref{eq:DiscreteBCCondition}.

\bibliographystyle{plain}
\bibliography{DGDG.bib}

\end{document}